\documentclass{article}
\usepackage{arxiv}
\usepackage[utf8]{inputenc} 
\usepackage[T1]{fontenc}    
\usepackage{hyperref}       
\usepackage{url}            
\usepackage{booktabs}       
\usepackage[english]{babel}
\usepackage{amsmath,amssymb,mathrsfs,stmaryrd,xcolor,theorem}
\usepackage{amsfonts}      
\usepackage{nicefrac}      
\usepackage{microtype}      
\usepackage{graphicx}
\usepackage{doi}
\usepackage{bm}
\usepackage[textsize=tiny,tickmarkheight=5pt]{todonotes}

\usepackage{geometry}
\usepackage{subcaption}
\usepackage{multirow}
\usepackage{comment}

\title{Energy-conserving scheme for a shallow-water model with capillary effects}

\author{
\href{https://orcid.org/0009-0005-1186-4167}{\includegraphics[scale=0.06]{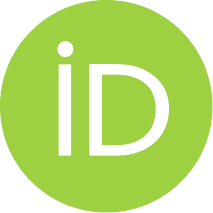}}\hspace{1mm}K. Sudarshan Kumar \\
	School of Mathematics\\
	Indian Institute of Science Education and Research\\
	Thiruvananthapuram-695551, India \\
	\texttt{sudarshan@iisertvm.ac.in} \\
	\And
	Boniface Nkonga \\
    Université Côte d’Azur\\
    INRIA/Castor, CNRS, LJAD\\
    06108 Nice Cedex 2, France\\
	\texttt{boniface.nkonga@univ-cotedazur.fr} \\
    \And
	\href{https://orcid.org/0000-0003-1903-4107}{\includegraphics[scale=0.06]{orcid.pdf}}\hspace{1mm}Praveen Chandrashekar\thanks{Corresponding author} \\
	Centre for Applicable Mathematics\\
    Tata Institute of Fundamental Research\\
    Bangalore -- 560065, India\\
	\texttt{praveen@tifrbng.res.in} \\
}

\hypersetup{
    colorlinks=true, 
    linkcolor=blue, 
    urlcolor=red, 
    citecolor=magenta,
    linktoc=all 
}

\newcommand{\assign}{:=}

\newcommand{\mathd}{\mathrm{d}}
\newcommand{\nobracket}{}

\newcommand{\tmcolor}[2]{{\color{#1}{#2}}}
\newcommand{\tmop}[1]{\ensuremath{\operatorname{#1}}}

\newcommand{\tmtextit}[1]{\text{{\itshape{#1}}}}

{\theorembodyfont{\rmfamily}\newtheorem{remark}{Remark}}

\newcommand{\od}[2]{\frac{\mathd #1}{\mathd #2}}

\newcommand{\jump}[1]{\llbracket #1 \rrbracket}

\newcommand{\avgc}[1]{\{\!\!\{#1\}\!\!\}}

\newcommand{\half}{\frac{1}{2}}

\newcommand{\jmh}{{j-\frac{1}{2}}}
\newcommand{\jph}{{j+\frac{1}{2}}}
\newcommand{\kmh}{{k-\frac{1}{2}}}
\newcommand{\kph}{{k+\frac{1}{2}}}
\newcommand{\nph}{{n+\frac{1}{2}}}

\newcommand{\Nmh}{{N-\frac{1}{2}}}

\newcommand{\BBb}[1]{\left [ #1 \right ]}
\newcommand{\BB}[1]{\left (  #1 \right )}

\def\E{{\cal E}}
\def\P{{\cal P}}
\def\H{{\cal H}}
\def\Ds{D}
\def\Us{U}
\def\Vs{V}
\def\Fs{F}
\def\Gs{G}
\def\Ms{M}

\def\TensS{A}
\def\TensV{\underline{\bm{A}}}

\def\gs{b}
\def\gV{\bm{\gs}}

\def\Ee{E}
\def\D{\bm{\Ds}}
\def\U{\bm{\Us}}
\def\V{\bm{\Vs}}
\def\F{\bm{\Fs}}
\def\M{\bm{\Ms}}
\def\Sm{\bm{S}}

\begin{document}
\maketitle
\begin{abstract}
We consider numerical approximations of shallow-water and thin-film flow models with capillary effects. These models lead to third-order differential equations, and we work with an augmented model with second-order derivatives. For a wide class of such models, which includes the classical shallow-water model and thin-film flow models under gravitational force, we develop a fully discrete energy-conserving scheme for both the quadratic and fully nonlinear capillary models. Numerical results show that the proposed scheme is stable and conserves total energy in the absence of dissipative effects.
\end{abstract}
\section{Introduction}

The shallow-water equations (SWE) provide a depth-averaged description of free-surface flows in the long-wave regime and were first derived in the context of open-channel hydraulics by Saint-Venant~\cite{SaintVenant1871}. Beyond geophysical and hydraulic applications, SWE-type models provide a useful setting for nonlinear waves and hyperbolic balance laws~\cite{Whitham1974,LeVeque2002}. Their numerical approximation has motivated a large literature on Riemann-solver-based finite-volume schemes and well-balanced treatments of source terms~\cite{Toro2009,Audusse2004,Fjordholm2011}, with extensions to additional physics such as magnetohydrodynamics~\cite{gilman2000,dellar2002} and surface tension~\cite{Serre1993,Casal1988}.

A recent application arises in magnetic confinement fusion, where a flowing liquid-metal blanket is being investigated as a means of shielding plasma-facing components from high heat fluxes~\cite{Abdou2015,Smolentsev2010}. The thin liquid-metal layer is subject to strong magnetic fields, and free-surface magnetohydrodynamic (MHD) flows of this type have been modeled within reduced shallow-water-type frameworks~\cite{Morley2004,Narula2006}. In this setting, capillary surface-tension effects are essential, yet few studies address the thin-film asymptotics for liquid metals; existing work includes experimental studies of open-channel liquid-metal flow~\cite{Fisher2018,Hvasta2018} and simple one-dimensional energy-balance analyses.

Here we develop numerical approximations of free-surface flows in the thin-film limit, accounting for capillary surface tension while neglecting the magnetic field. In the shallow-water approximation, long-wave scaling of the vertical momentum equation yields a balance between the vertical pressure gradient and gravity. Integrating this balance over the depth gives a pressure composed of a hydrostatic part and a free-surface contribution. The standard assumption that the pressure at the free surface equals the ambient pressure is valid only when surface-tension effects are negligible or the interface is flat. For liquid metals, capillary effects are not negligible, so the approximation must retain the pressure jump induced by capillary forces at the free surface. This jump is naturally incorporated through Hamilton's principle; the resulting Euler-Lagrange equations~\cite{Serre1993,Casal1988} contain third-order spatial derivatives. Discretizing third-order operators is challenging and often produces severe instabilities and spurious grid-scale oscillations~\cite{Bresch2020,Popinet2018,Favrie2017}; preserving the total-energy balance at the discrete level is equally difficult~\cite{Bresch2020}. Energy conservation is important in this setting because it prevents spurious energy growth in long-time integrations and allows dispersive capillary waves to be resolved without relying on artificial dissipation. A common remedy is to reformulate the model with lower-order derivatives. Bresch et al.~\cite{Bresch2020} introduce a scaled depth gradient as an auxiliary variable and obtain a second-order augmented system, while Dhaouadi et al.~\cite{Dhaouadi2022} use hyperbolic relaxation to obtain a first-order formulation.

The thin-film limit itself has an extensive literature. Benney~\cite{Benney1966} derived the classical long-wave film equation, and reduced models obtained by weighted residuals, such as the two-equation model of Ruyer-Quil and Manneville~\cite{RuyerQuil2000}, are now standard; see~\cite{Oron1997,Craster2009} for reviews and~\cite{Kalliadasis2012} for a monograph treatment.

In the formulation of~\cite{Bresch2020}, the scaled depth gradient makes the capillary energy quadratic and yields an augmented system that conserves total energy. Dhaouadi et al.~\cite{Dhaouadi2022} instead use an extended Lagrangian to express capillarity through lower-order derivatives, resulting in a first-order nonlinear system with a stiff source term and relaxation parameters that must be chosen by the user. We adopt the formulation of~\cite{Bresch2020} because its quadratic energy is amenable to entropy-stable discretization and its conservative form admits a discrete energy balance. The total energy plays the role of a convex entropy for the augmented system, which allows us to use the entropy-stability framework of Tadmor~\cite{Tadmor2003}. Our main contributions are:
\begin{itemize}
\item An energy-conserving two-point numerical flux, combined with summation-by-parts (SBP) finite differences~\cite{Mattsson2004,Svard2014}, which yields a provably energy-conserving semi-discrete scheme.
\item A Crank-Nicolson-type time discretization based on a special time-averaged state related to the system energy~\cite{Ray2017a}, which gives a fully discrete energy-conserving scheme.
\item An unsplit treatment of the hyperbolic and capillary terms. In contrast to the operator splitting used in~\cite{Bresch2020}, for which a rigorous discrete energy-conservation proof has not been established, the present method is provably energy-conserving.
\item A unified treatment of general hydrostatic potentials and of quadratic and fully nonlinear capillary energies in one and two dimensions.
\end{itemize}
The fully implicit nonlinear systems are solved by Newton's method with matrix-free generalized minimal residual (GMRES) for the linearized systems, using the Portable, Extensible Toolkit for Scientific Computation (PETSc). This nonlinear solver is the main computational cost of the method. It allows us to choose a time step based on the hyperbolic Courant-Friedrichs-Lewy (CFL) condition instead of the more restrictive condition imposed by dispersive effects.

We apply the method to two classes of problems: classical shallow water, with gravitational potential energy $\frac{1}{2} g h^2$, and thin-film flows~\cite{Richard2019}, whose potential energy has a more general polynomial form. Numerical experiments, including one- and two-dimensional Gaussian-bump problems and falling liquid films, demonstrate second-order convergence for the quadratic model and convergence rates approaching second order for the nonlinear model; they also show that the total energy is conserved to machine precision and that dispersive capillary waves are captured on sufficiently fine grids without added numerical dissipation. Since the numerical flux is of a central type, the method is intended for smooth solutions and may exhibit grid-scale oscillations on coarse grids, particularly for the nonlinear capillary model. Throughout, we assume a positive water depth, $h>0$, and either periodic or inflow/outflow boundaries. The present study concentrates on the capillary mechanism; coupling with the magnetic field and positivity preservation are left for future work.

The remainder of the paper is organized as follows. Section~\ref{sec:model} introduces the shallow-water model with capillary effects and its second-order augmented formulation. Section~\ref{sec:scheme1d} derives the energy-conserving scheme for the one-dimensional model, including its semi-discrete and fully discrete forms and its application to a two-equation thin-film model. Section~\ref{sec:scheme2d} extends the approach to the two-dimensional model. Numerical experiments are reported in Section~\ref{sec:numerics}, and Section~\ref{sec:conclusion} summarizes the main results and outlines future work. Appendix~\ref{sec:lgbcout} analyzes the outflow boundary condition used in the falling-film test case.
\section{Shallow-water model with capillary effects}\label{sec:model}
We denote by $h(\bm{x}, t) > 0$ the water depth and by $\bm{v}(\bm{x}, t) = (v_1, v_2)$ the two-dimensional depth-averaged velocity field.  Following the second-gradient theory~\cite{Casal1988}, the total internal energy $\E$ is expressed as a function of both the local depth $h$ and its spatial gradient $\nabla h$:
\begin{equation}\label{eq:iemodel}
\E(h, \nabla h) = \Phi(h) + \kappa(h) \E_{\tmop{cap}}(|\nabla h|)
\end{equation}
where $\Phi(h)$ accounts for the gravitational potential energy, and $\kappa(h)\E_{\tmop{cap}}(|\nabla h|)$ represents the extended chemical potential accounting for capillary effects. The capillary coefficient $\kappa(h) = \sigma(h)/\rho$ is positive, where $\sigma$ is the surface tension coefficient and $\rho$ is the fluid density. Applying the Euler-Lagrange equations under the variational constraint of mass conservation yields the shallow-water system with capillary surface tension~\cite{Bresch2020}:
\begin{equation}
  \begin{array}{rll}
    \partial_t h + \tmop{div} (h \bm{v}) & = & 0\\[2mm]
    \partial_t (h \bm{v}) + \tmop{div} (h \bm{v} \otimes
    \bm{v}) + \nabla \P& = & - \tmop{div} (\nabla h \otimes
    \nabla_{\bm{p}} \E) + \nabla (h \tmop{div} (\nabla_{\bm{p}} \E))
  \end{array} \label{eq:swc3}
\end{equation}
where the internal energy function $\E = \E (h, \bm{p})$ is given by \eqref{eq:iemodel} and
\[
\P(h, \nabla h) = h \partial_h \E (h, \nabla h) - \E (h, \nabla h) = \pi (h)
   - (\kappa (h) - h \kappa' (h)) \E_{\tmop{cap}} (| \nabla h |)
\]
The hydrostatic potential $\Phi(h)$ and the hydrostatic pressure component $\pi(h)$ are related by
\begin{equation}\label{eq:piphi}
\pi(h) = h\Phi'(h) - \Phi(h).
\end{equation}
To ensure that the system is physically admissible and mathematically well-posed, we enforce the following structural assumptions on the constitutive functions:
\begin{itemize}
\item The functions $s \mapsto \E_{\tmop{cap}}(s)$, $s \mapsto \Phi(s)$, and $s \mapsto \sigma(s)$ are strictly positive for $s > 0$.

\item The capillary energy function $\E_{\tmop{cap}}$ is invertible from $\mathbb{R}^+$ to $\mathbb{R}^+$ with $\E_{\tmop{cap}}(0) = 0$.

 \item The hydrostatic pressure satisfies $\pi'(h) = h\Phi''(h) > 0$, ensuring that $\Phi(h)$ is a strictly convex function for all admissible depths $h > 0$.


\end{itemize}
Two important models for capillarity are the {\em quadratic capillary} case for which the capillary energy function is given by
\begin{equation}
\E_{\tmop{cap}} (q) = \frac{1}{2} q^2 \label{eq:capquad}
\end{equation}
and the {\em fully nonlinear capillary} case, for which
\begin{equation}
\E_{\tmop{cap}} (q) = \sqrt{1 + q^2} - 1 \label{eq:capnonl}
\end{equation}
The quadratic model is an approximation of the fully nonlinear model when $q \approx 0$. The model \eqref{eq:swc3} admits an energy conservation equation
\[
\partial_t \left( \frac{1}{2} h | \bm{v} |^2 + \E \right) + \tmop{div}
   \left[ \left( \frac{1}{2} h | \bm{v} |^2 + \E + \P\right) \bm{v}
- \bm{Q} \right] = 0
\]
where $ \bm{Q} = \tmop{div} (\nabla_{\bm{p}} \E) h \bm{v}
- \tmop{div} (h \bm{v}) \nabla_{\bm{p}} \E $.
\paragraph{Reformulation as a second-order system.}
Bresch et al.~\cite{Bresch2020} treat the depth gradient as an independent auxiliary vector-field variable and define
\[ \bm{w} = \BBb{\begin{array}{c} w_1\\ w_2\end{array} }= \alpha (| \nabla h |^2) \sqrt{\frac{\kappa
   (h)}{h}} \nabla h, \qquad \alpha (s) =
   \sqrt{\frac{2\E_{\tmop{cap}} \left( \sqrt{s} \right)}{s}}
\]
This new variable, which has units of velocity, allows us to recast the capillary internal energy as a quadratic term which is structurally analogous to classical kinetic energy, i.e.,
\[
\frac{1}{2} h | \bm{w} |^2 = \kappa (h) \E_{\tmop{cap}} (| \nabla h |) = \text{capillary energy}
\]
Using the continuity equation it is possible to derive a partial differential equation for $\bm{w}$. The resulting second-order augmented hyperbolic system takes the form
\begin{equation}
  \U_t + \tmop{div} \F (\U) = \M(\U), \qquad \U = \left[\begin{array}{c}
    h\\
    h \bm{v}\\
    h \bm{w}
  \end{array}\right], \qquad \F = \left[\begin{array}{c}
    h \bm{v}\\
    h \bm{v} \otimes \bm{v} + \pi (h) I\\
    h \bm{v} \otimes \bm{w}
  \end{array}\right], \qquad \M = \left[\begin{array}{c}
    0\\
    \M_v\\
    \M_w
  \end{array}\right] \label{eq:aug2d}
\end{equation}
where
\begin{equation}
  \M_v = \tmop{div} \BB{h \nabla (\TensV\bm{w})^{\top}} - \nabla
  (\bm{\gs} \cdot \bm{w}), \qquad \M_w = - \TensV \tmop{div} \BB{h
  \nabla \bm{v}^{\top}} - \bm{\gs} \tmop{div} \bm{v}
  \label{eq:M23}
\end{equation}
with $\TensV \equiv \TensV(h, \bm{w})$ a symmetric tensor and $\bm{\gs} \equiv \bm{\gs} (h, \bm{w})$ a vector field. This model satisfies the energy/entropy equation
\[
\partial_t \Ee + \tmop{div} [(\Ee + \pi) \bm{v} - \bm{R}] = 0 \qquad \text{ where } \qquad \Ee = \frac{1}{2} h |\bm{v}|^2 + \frac{1}{2} h
   |\bm{w}|^2 + \Phi (h)
\]
where $\bm{R} = \nabla (\TensV\bm{w}) h \bm{v} - h (\nabla\bm{v}) \TensV \bm{w} - (\bm{\gs} \cdot \bm{w}) \bm{v}$. In the following sections, we first develop the energy-conserving scheme for the one-dimensional model and then extend it to the two-dimensional case.

\section{Energy-conserving scheme for the augmented 1-D model}\label{sec:scheme1d}
We consider the one-dimensional version of the augmented model~\eqref{eq:aug2d},~\eqref{eq:M23} and give a detailed derivation of the numerical scheme. The quantities $\bm{v}, \bm{w}$ have a single component each, denoted by scalar functions $v, w$ respectively. The augmented model is then given by
\begin{equation} \label{eq:bresch1d}
  \U_t + \partial_x \F (\U) = \M (\U)
\end{equation}
where
\[ \U = \left[\begin{array}{c}
     h\\[1mm]
     h v\\[1mm]
     h w
   \end{array}\right], \qquad \F = \left[\begin{array}{c}
     h v\\[1mm]
     h v^2 + \pi (h)\\[1mm]
     h v w
   \end{array}\right], \qquad \M = \left[\begin{array}{c}
     0\\[1mm]
     \partial_x (h \partial_x (\TensS w)) - \partial_x (\gs w)\\[1mm]
     - \TensS \partial_x (h \partial_x v) - \gs \partial_x v
   \end{array}\right] = \left[\begin{array}{c}
     0\\[1mm]
     \Ms_v\\[1mm]
     \Ms_w
   \end{array}\right]
   \]
and the scalar functions $\TensS$ and $\gs$ are defined as
\begin{eqnarray*}
  \TensS(h, w) & = & \sqrt{h \kappa (h)} \left[ 2 \frac{\alpha' (q^2) h}{\alpha^2
  (q^2) \kappa (h)} w^2 + \alpha (q^2) \right]\\
  \gs(h, w) & = & \left[ \frac{\kappa' (h) h}{2 \kappa (h)} + \frac{1}{2} + 2
  \frac{\alpha' (q^2)}{\alpha (q^2)} q^2 \right] h w\\
    q & = & \E^{- 1}_{\tmop{cap}} \left( \frac{h w^2}{2 \kappa (h)}
  \right)
\end{eqnarray*}
The flux Jacobian is
\[ \F' (\U) = \left[\begin{array}{ccc}
     0 & 1 & 0\\
     \pi' - v^2 & 2 v & 0\\
     - v w & w & v
   \end{array}\right] \]
The eigenvalues of the flux Jacobian are $v - \sqrt{\pi' (h)}, v, v +
\sqrt{\pi' (h)}$ and they are real by convexity of $\Phi (h)$ for $h > 0$. The
corresponding eigenvectors
\[ \left[\begin{array}{c}
     1\\
     v - \sqrt{\pi' (h)}\\
     w
   \end{array}\right], \qquad \left[\begin{array}{c}
     0\\
     0\\
     1
   \end{array}\right], \qquad \left[\begin{array}{c}
     1\\
     v + \sqrt{\pi' (h)}\\
     w
   \end{array}\right] \]
are linearly independent. The model admits the total energy
\[ \Ee = \frac{1}{2} h v^2 + \frac{1}{2} h w^2 + \Phi (h) \]
and satisfies the energy/entropy equation~\cite{Bresch2020}
\[ \partial_t \Ee + \partial_x\Gs = \partial_x\BB{h v \partial_x (\TensS w) - \TensS w h \partial_x v - \gs v w},
   \qquad \Gs = (\Ee + \pi) v \]
as verified by a direct calculation.

An alternative approach to obtain this energy equation is to use the framework of entropy function and entropy variables. Following the well-established idea of Tadmor~\cite{Tadmor2003}, let us introduce the \tmtextit{entropy variables}
\[ \V = \Ee' (\U)^{\top} = \left[\begin{array}{c}
     \Phi' (h) - \frac{1}{2} (v^2 + w^2)\\
     v\\
     w
   \end{array}\right] \]
Moreover,
\begin{eqnarray*}
  \V^{\top} \F' (\U) & = & \left[\begin{array}{ccc}
    \Phi' (h) - \frac{1}{2} (v^2 + w^2) & v & w
  \end{array}\right] \left[\begin{array}{ccc}
    0 & 1 & 0\\
    \pi' (h) - v^2 & 2 v & 0\\
    - v w & w & v
  \end{array}\right]\\
  & = & \left[\begin{array}{ccc}
    v \pi' (h) - v^3 - v w^2, & \Phi' (h) + \frac{3}{2} v^2 + \frac{1}{2} w^2,
    & v w
  \end{array}\right]\\
  & = & \Gs' (\U)
\end{eqnarray*}
which shows that $(\Ee (\U), \Gs (\U))$ is a {\em convex entropy pair} satisfying
\[ \Ee' (\U) \F' (\U) = \Gs' (\U) \]
We can thus obtain the energy equation by taking the dot product of \eqref{eq:bresch1d} with $\V$;
\[ \V \cdot (\partial_t\U + \partial_x\F) = \V \cdot \M \qquad \Longrightarrow \qquad \partial_t \Ee + \partial_x\Gs = \V \cdot \M \]
Moreover,
\begin{eqnarray*}
  \V \cdot \M & = & \tmcolor{blue}{v \partial_x (h \partial_x (\TensS w))}
  \tmcolor{red}{- v \partial_x (\gs w)} - \TensS w \partial_x (h \partial_x v)
  \tmcolor{red}{- \gs w \partial_x v}\\
  & = & \tmcolor{blue}{\partial_x (h v \partial_x (\TensS w)) - h \partial_x (\TensS w)
  \partial_x v} - \partial_x (\TensS w h \partial_x v) + \partial_x (\TensS w) h
  \partial_x v \tmcolor{red}{- \partial_x (\gs v w)}\\
  & = & \partial_x (h v \partial_x (\TensS w)) - \partial_x (\TensS w h \partial_x v) -
  \partial_x (\gs v w)\\
  & = & \partial_x [h v \partial_x (\TensS w) - \TensS w h \partial_x v - \gs v w]
\end{eqnarray*}
and the last expression on the right is the divergence of some flux. Under periodic boundary conditions or assuming no net boundary flux, we obtain
\[ \int_{\Omega} \partial_t\Ee \mathd x = - \int_{\Omega} \partial_x\Gs \mathd x + \int_{\Omega} \V
   \cdot \M \mathd x = 0 \qquad \Longrightarrow \qquad \int_{\Omega} \Ee (x, t)
   \mathd x = \tmop{const} \]
Thus the total energy in the domain is constant with respect to time.

To obtain this property in the finite volume method, we have to construct a numerical flux to satisfy certain condition. Define the \tmtextit{entropy potential} by
\begin{eqnarray*}
  \H & = & \V \cdot \F - \Gs\\
  & = & \left[\begin{array}{ccc}
    \Phi' (h) - \frac{1}{2} (v^2 + w^2) & v & w
  \end{array}\right] \left[\begin{array}{c}
    h v\\
    h v^2 + \pi (h)\\
    h v w
  \end{array}\right] - (\Ee + \pi) v\\
  & = & [h \Phi' (h) - \Phi (h)] v\\
  & = & \pi (h) v
\end{eqnarray*}
Assume that there is a two-point numerical flux $\F^{\star} (\U_l, \U_r)$ which satisfies the following condition~\cite{Tadmor1984,Tadmor2003}
\begin{equation}
  (\V_r - \V_l) \cdot \F^{\star} (\U_l, \U_r) = \H (\U_r) - \H (\U_l) \label{eq:Fsdef}
\end{equation}
for any two physically admissible states $\U_l, \U_r$ with $\V_l = \V(\U_l)$ and $\V_r = \V(\U_r)$.  This is a single equation for the three components of $\F^{\star}$ and we will discuss the form of such a flux in the later sections and construct formulae for their numerical implementation. 
\subsection{Semi-discrete scheme: hyperbolic part}\label{sec:sdhyp}
The semi-discrete scheme using the flux $\F^{\star}$ is given by
\[
\od{\U_j}{t} + \frac{\F^{\star} (\U_j, \U_{j + 1}) - \F^{\star} (\U_{j - 1}, \U_j)}{\Delta x} = 0
\]
where we have ignored the second-order terms $\M$; we will add them back in a later section. Taking the dot product of entropy variables $\V_j$ with the above equation, we obtain the semi-discrete energy equation
\[
\od{\Ee (\U_j)}{t} + \frac{\Gs^{\star} (\U_j, \U_{j + 1}) - \Gs^{\star} (\U_{j - 1},
   \U_j)}{\Delta x} = 0
\]
with a {\em numerical energy flux} $\Gs^{\star}$ given by
\begin{equation}\label{eq:Fes}
\Gs^{\star} (\U_j, \U_{j+1}) = \left( \frac{\V_j + \V_{j + 1}}{2} \right) \cdot
   \F^{\star} (\U_j, \U_{j + 1}) - \frac{\H_j + \H_{j + 1}}{2}
\end{equation}
Clearly, under periodic boundary conditions, the total energy in the spatial domain, $\Delta x \sum_j \Ee (\U_j)$, is constant in time.

\subsection{Energy-conserving numerical flux}
To derive an energy-conserving flux, we expand the discrete entropy-conserving condition~\eqref{eq:Fsdef} in terms of an independent set of variables and solve for the unknown flux components using a system of equations for each of the jump terms. Define the {\em average} and {\em jump} operators corresponding to any two states which are indexed by $l, r$ as
\[ \avgc{\cdot} = \frac{1}{2} [(\cdot)_l + (\cdot)_r], \qquad \jump{\cdot} =
   (\cdot)_r - (\cdot)_l \]
respectively. We will make repeated use of the following elementary relation for the jump of a product of two quantities $a,b$,
\[
\jump{ab} = \avgc{a}\jump{b} + \avgc{b}\jump{a}
\]
\paragraph{$\F^{\star}$ when $\Phi (h) = \frac{1}{2} g h^2$.} This corresponds to the usual shallow-water model and in this case $\pi (h) = \frac{1}{2} g h^2$ and the entropy variables are given by
\[
\V = \Ee' (\U)^{\top} = \left[\begin{array}{c}
     g h - \frac{1}{2} (v^2 + w^2)\\
     v\\
     w
   \end{array}\right]
\]
and the entropy potential is
\[
\H = \pi (h) v = \frac{1}{2} g h^2 v
\]
The entropy conservation condition \eqref{eq:Fsdef} for the flux $\F^{\star}$ yields the condition,
\[
\left[\begin{array}{c}
     g \jump{h} - \avgc{v} \jump{v} - \avgc{w} \jump{w}\\ \jump{v}\\
     \jump{w}
   \end{array}\right] \cdot \left[\begin{array}{c}
     \Fs_1^{\star}\\
     \Fs_2^{\star}\\
     \Fs_3^{\star}
   \end{array}\right] = \jump{\H} = g \avgc{h} \avgc{v} \jump{h} + \frac{1}{2}
   g \avgc{h^2} \jump{v}
\]
i.e.,
\[ g \BB{\rule{0mm}{4mm} \Fs_1^{\star} - \avgc{h} \avgc{v}} \jump{h}
   - \BB{ \rule{0mm}{4mm} \avgc{v}\Fs_1^{\star} -  \Fs_2^{\star} - \frac{g}{2}\avgc{h^2}} \jump{v}
   - \BB{\rule{0mm}{4mm}\avgc{w} \Fs_1^{\star}
   - \Fs_3^{\star} } \jump{w}  = 0
\]
Since $h, v, w$ are independent variables, the above equation must hold for all values of $\jump{h}, \jump{v}, \jump{w}$; thus the coefficients of each of the jump terms must vanish, which uniquely determines all the components of the flux,
\[
\F^{\star} \equiv \F^{\star}\BB{\U_l, \U_r} = \left[\begin{array}{c}
     \avgc{h} \avgc{v}\\[2mm]
     \avgc{v} \Fs_1^{\star} + \frac{1}{2} g \avgc{h^2}\\[2mm]
     \avgc{w} \Fs_1^{\star}
   \end{array}\right] = \left[\begin{array}{c}
     \avgc{h} \avgc{v}\\[2mm]
     \avgc{h} \avgc{v}^2 + \frac{1}{2} g \avgc{h^2}\\[2mm]
     \avgc{h} \avgc{v} \avgc{w}
   \end{array}\right]
\]
This numerical flux is consistent since $\F^{\star} (\U, \U) = \F (\U)$ and is a smooth function of its arguments.

\paragraph{$\F^{\star}$ for general $\Phi (h)$.}

The entropy variables and entropy potential are given by
\begin{equation}\label{eq:VHgenPhi}
\V = \Ee' (\U)^{\top} = \left[\begin{array}{c}
     \Phi' (h) - \frac{1}{2} (v^2 + w^2)\\
     v\\
     w
   \end{array}\right] \qquad \tmop{and} \qquad \H = \pi(h) v
\end{equation}
respectively. The condition \eqref{eq:Fsdef} for entropy-conserving flux
becomes
\[
\left[\rule{0mm}{4mm}
\begin{array}{c}
     \jump{\Phi'} - \avgc{v} \jump{v} - \avgc{w} \jump{w}\\  \jump{v}\\
     \jump{w}
   \end{array}\right] \cdot \left[\begin{array}{c}
     \Fs_1^{\star}\\
     \Fs_2^{\star}\\
     \Fs_3^{\star}
   \end{array}\right] = \jump{\H} = \avgc{v} \jump{\pi} + \avgc{\pi} \jump{v}
\]
Since $\Phi', \pi$ depend only on $h$, we can write their jumps in terms of
the jump in $h$, and the above equation becomes
\begin{equation*}
\BB{\rule{0mm}{4mm}{\frac{\jump{\Phi'}}{\jump{h}} \Fs_1^{\star}  - \avgc{v} \frac{\jump{\pi}}{\jump{h}}}} \jump{h} + \left(
   \rule{0mm}{4mm}{- \avgc{\pi} - \avgc{v} \Fs_1^{\star} + \Fs_2^{\star}} \right) \jump{v} +
   \left( \rule{0mm}{4mm}- \avgc{w} \Fs_1^{\star} + \Fs_3^{\star} \right) \jump{w} = 0
\end{equation*}
This equality should hold for any values of the jumps $\jump{h}$, $\jump{v} $ and $\jump{w}$.  Therefore, the flux component should be defined as
\begin{equation}
  \F^{\star} =
    \left[\begin{array}{c}
    \Fs_1^{\star} \\[3mm]
   \Fs_2^{\star} \\[3mm]
    \Fs_3^{\star}
  \end{array}\right]  =
  \left[\begin{array}{c}
    \frac{\jump{\pi}}{\jump{\Phi'}} \avgc{v}\\[3mm]
    \avgc{v} \Fs_1^{\star} + \avgc{\pi}\\[3mm]
    \avgc{w} \Fs_1^{\star}
  \end{array}\right]  =
  \left[\begin{array}{c}
    \hat{h} \avgc{v}\\[3mm]
    \hat{h} \avgc{v}^2  + \avgc{\pi}\\[3mm]
    \hat{h} \avgc{w} \avgc{v}
  \end{array}\right], \qquad \hat h \assign \frac{\jump{\pi}}{\jump{\Phi'}}
  \label{eq:Fstfm}
\end{equation}
Note that since $\pi(h) = h \Phi' (h) - \Phi (h)$, we have
\begin{equation}
  \hat{h} = \avgc{h} +
  \frac{\avgc{\Phi'} - \jump{\Phi} / \jump{h}}{\jump{\Phi'} / \jump{h}}
  \label{eq:hh}
\end{equation}
If $h_l = h_r = h$, then
\[ \hat{h} = h + \frac{\Phi' (h) - \Phi' (h)}{\Phi'' (h)} = h \qquad
   \tmop{since} \qquad \Phi'' > 0 \]
Thus $\F^{\star}$ is a consistent numerical flux. Note that $\pi = \pi (h)$,
$\Phi = \Phi (h)$ and we will assume that for all the models of interest, the
quantities $\jump{\pi} / \jump{h}$, \ $\jump{\Phi} / \jump{h}$ and $\jump{\Phi'}
/ \jump{h}$ are well defined for any two admissible states $h_l, h_r$. This
also implies that $\hat{h} = \frac{\jump{\pi}}{\jump{\Phi'}}$ is also well
defined. Such a model of physical interest for flow of a thin film down an
inclined plane is discussed in Section~\ref{sec:richard}.

\subsection{Fully discrete scheme: hyperbolic part}\label{sec:fdhyp}
Consider a time average
\[
\bar{\U} = \bar{\U} (\U^{(1)}, \U^{(2)})
\]
where $\U^{(1)}, \U^{(2)}$ are the state variables at two different times. We will assume that the above time average has the property that~\cite{Ray2017a}
\begin{equation}
  \V \BB{\bar{\U} \BB{\U^{(1)}, \U^{(2)}}} \cdot \BB{\U^{(2)} - \U^{(1)}} = \Ee \BB{\U^{(2)}}- \Ee\BB{\U^{(1)}} \label{eq:tavg}
\end{equation}
and consider the Crank-Nicolson-type scheme
\[ \frac{\U^{n + 1}_j - \U_j^n}{\Delta t} + \frac{\F^{\star}\BB{\U_j^{\nph},
   \U_{j + 1}^{\nph}} - \F^{\star}\BB{\U_{j - 1}^{\nph}, \U_j^{\nph}}}{\Delta x} = 0, \qquad \U^{\nph}_j = \bar{\U} \BB{\U^n_j, \U^{n + 1}_j}
\]
where we have again ignored the second-order terms in $\M$. Taking the dot product with $\V (\U^{\nph}_j)$ and using~\eqref{eq:tavg}, we get
\[ \frac{\Ee (\U^{n + 1}_j) - \Ee (\U^n_j)}{\Delta t}
+ \frac{  \Gs^{\star} \BB{\U_j^{\nph}, \U_{j + 1}^{\nph}}
        - \Gs^{\star} \BB{\U_{j - 1}^{\nph}, \U_j^{\nph}}
        }{\Delta x} = 0
        \]
where the numerical energy flux $\Gs^\star$ is given by~\eqref{eq:Fes}. Thus the fully discrete scheme conserves total energy under periodic boundary
conditions.

\paragraph{$\bar{\U}$ when $\Phi (h) = \frac{1}{2} g h^2$.}

In the following, $\avgc{\cdot}$ and $\jump{\cdot}$ denote averages and jumps in the time direction; the subscripts $l, r$ now denote time levels $n, n + 1$. We have
\[ \jump{\U} = \left[\begin{array}{c}
     \jump{h}\\[2mm]
     \avgc{v} \jump{h} + \avgc{h} \jump{v}\\[2mm]
     \avgc{w} \jump{h} + \avgc{h} \jump{w}
   \end{array}\right] \]
and
\[ \jump{\Ee} = \left( \tmcolor{blue}{g \avgc{h} + \frac{1}{2} \avgc{v^2} +
   \frac{1}{2} \avgc{w^2}} \right) \jump{h} + \tmcolor{red}{\avgc{h} \avgc{v}}
   \jump{v} + \tmcolor{magenta}{\avgc{h} \avgc{w}} \jump{w} \]
Then, with $\bar{\V} = \V (\bar{\U} (\U_l, \U_r))$, the condition \eqref{eq:tavg}
becomes $\bar{\V} \cdot \jump{\U} = \jump{\Ee}$ and yields
\[ \left( \tmcolor{blue}{\bar{\Vs}_1 + \avgc{v} \bar{\Vs}_2 + \avgc{w} \bar{\Vs}_3}
   \right) \jump{h} + \tmcolor{red}{\avgc{h} \bar{\Vs}_2} \jump{v} +
   \tmcolor{magenta}{\avgc{h} \bar{\Vs}_3} \jump{w} = \jump{\Ee} \]
Comparing the corresponding jump terms on both sides, we get
\begin{eqnarray*}
  \bar{\Vs}_2 & = & \avgc{v}\\[2mm]
  \bar{\Vs}_3 & = & \avgc{w}\\
  \bar{\Vs}_1 & = & g \avgc{h} + \frac{1}{2} \avgc{v^2} + \frac{1}{2} \avgc{w^2}
  - \bar{\Vs}_2^2 - \bar{\Vs}_3^2
   \;\; = \;\;  g \avgc{h} - \frac{1}{2} v_l v_r - \frac{1}{2} w_l w_r
\end{eqnarray*}
We have to convert the entropy variables $\bar{\V}$ to conserved ($\bar{\U}$) or
primitive variables $(\bar h, \bar v, \bar w)$. The map $\V \rightarrow (h, v, w)$ is given by
\[ h = \frac{1}{g} \left( \Vs_1 + \frac{\Vs_2^2}{2} + \frac{\Vs_3^2}{2} \right),
   \qquad v = \Vs_2, \qquad w = \Vs_3 \]
In fact, $\bar{v} = \avgc{v}$, $\bar{w} = \avgc{w}$ and
\[ \bar{h} = \frac{1}{g} \left( \bar{\Vs}_1 + \frac{\bar{\Vs}_2^2}{2} +
   \frac{\bar{\Vs}_3^2}{2} \right) = \avgc{h} + \frac{1}{8 g} [(v_r - v_l)^2 +
   (w_r - w_l)^2] > 0 \]
provided $h_l, h_r > 0$. Thus the time average is well defined and consistent
if the two states are physically admissible.

\paragraph{$\bar{\U}$ for general $\Phi (h)$.}

For a general potential energy function $\Phi (h)$, the jump in total energy
can be written as
\[ \jump{\Ee} = \left( \tmcolor{blue}{\frac{\jump{\Phi}}{\jump{h}} + \frac{1}{2}
   \avgc{v^2} + \frac{1}{2} \avgc{w^2}} \right) \jump{h} +
   \tmcolor{red}{\avgc{h} \avgc{v}} \jump{v} + \tmcolor{magenta}{\avgc{h}
   \avgc{w}} \jump{w} \]
Then, with $\bar{\V} = \V\BB{\bar{\U} (\U_l, \U_r)}$, the condition \eqref{eq:tavg}
becomes $\bar{\V} \cdot \jump{\U} = \jump{\Ee}$ and yields
\[ \left( \tmcolor{blue}{\bar{\Vs}_1 + \avgc{v} \bar{\Vs}_2 + \avgc{w} \bar{\Vs}_3}
   \right) \jump{h} + \tmcolor{red}{\avgc{h} \bar{\Vs}_2} \jump{v} +
   \tmcolor{magenta}{\avgc{h} \bar{\Vs}_3} \jump{w} = \jump{\Ee} \]
and comparing jump terms, we get
\begin{eqnarray*}
  \bar{\Vs}_2 & = & \avgc{v}\\[2mm]
  \bar{\Vs}_3 & = & \avgc{w}\\[2mm]
  \bar{\Vs}_1 & = & \frac{\jump{\Phi}}{\jump{h}} + \frac{1}{2} \avgc{v^2} +
  \frac{1}{2} \avgc{w^2} - \bar{\Vs}_2^2 - \bar{\Vs}_3^2 \;\; = \;\; \frac{\jump{\Phi}}{\jump{h}} - \frac{ v_l v_r}{2} - \frac{w_lw_r}{2}
\end{eqnarray*}
We have to convert the entropy variables $\bar{\V}$ to conserved variables $\bar{\U}$ or primitive variables $(\bar h, \bar v, \bar w)$. The map $\V \rightarrow (h, v, w)$ is given by
\[ h = (\Phi')^{- 1} \left( \Vs_1 + \frac{\Vs_2^2}{2} + \frac{\Vs_3^2}{2} \right),
   \qquad v = \Vs_2, \qquad w = \Vs_3 \]
In fact, from \eqref{eq:VHgenPhi}
\begin{eqnarray}
  \Phi' (\bar{h}) & = & \bar{\Vs}_1 + \frac{\bar{\Vs}_2^2}{2} +
  \frac{\bar{\Vs}_3^2}{2}  \;\; = \;\;  \frac{\jump{\Phi}}{\jump{h}} - \frac{1}{2} v_l v_r - \frac{1}{2} w_l
  w_r + \frac{1}{2} \left( \frac{v_l + v_r}{2} \right)^2 + \frac{1}{2} \left(
  \frac{w_l + w_r}{2} \right)^2 \nonumber\\
  & = & \frac{\jump{\Phi}}{\jump{h}} + \frac{(v_l - v_r)^2 }{8} +
  \frac{(w_l - w_r)^2}{8}   \label{eq:barh}
\end{eqnarray}
which is an implicit equation for $\bar{h}$. Now, since $\Phi$ is convex,
$\Phi'$ is an increasing function for $h > 0$; in fact it is reasonable to
assume that $\Phi' (0) = 0$. If $h_l, h_r > 0$ then $\jump{\Phi} / \jump{h} >
0$, and the right hand side of \eqref{eq:barh} is positive, so that there is a
unique solution $\bar{h} > 0$ for any two states $h_l, h_r > 0$. In the
computer code, we can use Newton's method to compute the average state
$\bar{h}$.

\subsection{Discretization of capillary terms}\label{sec:Mdisc}

So far we have discussed how the hyperbolic part of the model \eqref{eq:bresch1d} will be approximated so that the total energy is conserved. We next discretize the second-order terms $\M$ by using standard central differencing, leading to
\begin{eqnarray*}
  (\Ms_v)_j & = & \frac{1}{\Delta x} \BB{ h_{\jph} \frac{(\TensS w)_{j + 1} -
  (\TensS w)_j}{\Delta x} - h_{\jmh} \frac{(\TensS w)_j - (\TensS w)_{j - 1}}{\Delta x} }
  - \frac{(\gs w)_{j + 1} - (\gs w)_{j - 1}}{2 \Delta x}\\[3mm]
  (\Ms_w)_j & = & - \frac{\TensS_j}{\Delta x} \BB{ h_{\jph} \frac{v_{j + 1} -
  v_j}{\Delta x} - h_{\jmh} \frac{v_j - v_{j - 1}}{\Delta x} } -
  \gs_j \frac{v_{j + 1} - v_{j - 1}}{2 \Delta x}
\end{eqnarray*}
where
\[ h_{\jph} = \frac{h_j + h_{j + 1}}{2} \]
To examine how these terms contribute to the energy equation, we use some standard results on summation-by-parts property for the second-order terms. Define
\[
\delta a_j = \frac{a_{\jph} - a_{\jmh}}{\Delta x}, \qquad \mu
   a_{\jph} = \frac{a_j + a_{j + 1}}{2}
\]
Then we can write
\[
(\Ms_v)_j = \delta ( (\mu h) \delta (\TensS w))_j - \delta \mu (\gs w)_j, \qquad
   (\Ms_w)_j = - \TensS_j \delta((\mu h) \delta v)_j - \gs_j \delta \BB{\mu v_j}
\]
We will assume periodic boundary conditions.  Now
\begin{eqnarray*}
  \sum_j \V_j \cdot \M_j & = &
  \sum_j \left ( \rule{0mm}{4mm} v_j  (\Ms_v)_j + w_j  (\Ms_w)_j\right )\\
  & = & \sum_j \left ( \rule{0mm}{4mm}\tmcolor{blue}{v_j \delta [(\mu h) \delta (\TensS w)]_j}
  \tmcolor{red}{- v_j \delta \mu (\gs w)_j} \tmcolor{blue}{- \TensS_j w_j \delta
  [(\mu h) \delta v]_j} \tmcolor{red}{- \gs_j w_j \delta \mu v_j} \right )
  \;\; = \;\;  0
\end{eqnarray*}
where we used \eqref{eq:sbp1} for the red terms and \eqref{eq:sbp2} for the blue terms. Thus the second order terms in $\M$ which arise due to capillary effects do not change the total energy in the domain which is consistent with what happens at the continuous level.

\subsection{Fully discrete unsplit scheme}
For the model \eqref{eq:bresch1d} and now including both hyperbolic and capillary terms, consider the Crank-Nicolson scheme
\begin{equation}
  \frac{\U^{n + 1}_j - \U_j^n}{\Delta t} + \frac{\F^{\star}\BB{\U_j^{\nph},
  \U_{j + 1}^{\nph}} - \F^{\star}\BB{ \U_{j - 1}^{\nph}, \U_j^{\nph}}}{\Delta x} = \M_j \BB{\U^{\nph}} \label{eq:fdcn}
\end{equation}
with the time average $\U^{\nph}_j = \bar{\U}\BB{\U^n_j, \U^{n + 1}_j}$ satisfying \eqref{eq:tavg} and numerical flux satisfying \eqref{eq:Fsdef}. Taking the dot product of the above scheme with $\V\BB{\U^{\nph}_j}$, summing over all cells and using the results from Sections~\ref{sec:sdhyp},~\ref{sec:fdhyp},~\ref{sec:Mdisc}, we obtain the energy conservation property
\[
\Delta x \sum_j \Ee (\U_j^{n + 1}) = \Delta x \sum_j \Ee (\U_j^n)
\]
under periodic boundary conditions.

\subsection{Application to a two-equation model for thin-film flows}
\label{sec:richard}
Richard et al.~\cite{Richard2019} derived a model for thin-film flows down
an inclined plane, which can be cast into the framework of the Bresch et al.
model~\cite{Bresch2020}. In dimensional units, the equations are given by
\begin{eqnarray*}
  \partial_t h + \partial_x (h v) & = & 0\\
  \partial_t (h v) + \partial_x (h v^2 + \pi (h)) & = & g \sin (\theta) h -
  \frac{3 \nu v}{h} + \frac{9 \nu}{2} \partial_x (h \partial_x v) +
  \frac{\sigma}{\rho} h \partial_x K
\end{eqnarray*}
where
\begin{equation}
  \pi (h) = \frac{2 g^2 \sin^2 \theta}{225 \nu^2} h^5 + \frac{g \cos
  (\theta)}{2} h^2, \qquad \theta \in [0, \pi / 2] \label{eq:richpi}
\end{equation}
From \eqref{eq:piphi}, this corresponds to the potential energy function
\begin{equation}
  \Phi (h) = \frac{g^2 \sin^2 \theta}{450 \nu^2} h^5 + \frac{g \cos
  (\theta)}{2} h^2 \label{eq:richPhi}
\end{equation}
and
\[ K = \partial_{x x} h \qquad \tmop{or} \qquad K = \partial_x \left(
   \frac{\partial_x h}{\sqrt{1 + (\partial_x h)^2}} \right) \]
corresponding to quadratic \eqref{eq:capquad} and nonlinear capillarity
\eqref{eq:capnonl}, respectively. This model can be cast in the form of
\eqref{eq:bresch1d} leading to some extra terms $\Sm$ in the momentum equation,
\[ \U_t + \F (\U)_x = \M (\U) + \Sm (\U), \qquad \Sm (\U) = \left[\begin{array}{c}
     0\\
     g \sin (\theta) h - \frac{3 \nu v}{h} + \frac{9 \nu}{2} \partial_x (h
     \partial_x v)\\
     0
   \end{array}\right] \]
The fully discrete Crank-Nicolson scheme for the above model is similar to
\eqref{eq:fdcn} and is given by
\[ \frac{\U^{n + 1}_j - \U_j^n}{\Delta t} + \frac{\F^{\star} (\U_j^{\nph},
   \U_{j + 1}^{\nph}) - \F^{\star} (\U_{j - 1}^{\nph}, \U_j^{\nph})}{\Delta x} = \M_j (\U^{\nph}) + \Sm_j (\U^{\nph}) \]
and the extra terms are discretized using a central difference scheme
\[ (\Sm_j (\U))_2 = g \sin (\theta) h_j - \frac{3 \nu v_j}{h_j} + \frac{9 \nu}{2}
   \BB{ \frac{h_{\jph} \frac{v_{j + 1} - v_j}{\Delta x} - h_{\jmh}
   \frac{v_j - v_{j - 1}}{\Delta x}}{\Delta x} }\]
To derive the energy equation, we take the dot product with $\V_j = \V (\U^{\nph}_j)$ and sum over all cells, leading to
\begin{eqnarray*}
  \Delta x \sum_j \Ee \BB{\U_j^{n + 1}} & = & \Delta x \sum_j \Ee \BB{\U_j^n} + g \sin
  (\theta) \Delta t \sum_j h_j^{\nph} v_j^{\nph} \Delta x\\
  &  & - 3 \nu \Delta t \sum_j \frac{(v_j^\nph)^2}{h_j^\nph} \Delta x - \frac{9 \nu\Delta t}{2} \sum_j
  h_{\jph}^{\nph} \left( \frac{v_{j + 1}^{\nph} - v_j^{\nph}}{\Delta x} \right)^2 \Delta x
\end{eqnarray*}
under periodic boundary conditions. The last two terms on the right, which add friction and viscosity to the model, have a negative sign and lead to loss of total energy under the assumption that $h$ remains positive. Thus the proposed numerical scheme is consistent with the energy equation at the continuous level.

\paragraph{Computation of space and time averages.} In the two-point flux \eqref{eq:Fstfm} for the thin-film model, we need the spatial average $\hat{h}$ which can be computed from \eqref{eq:hh} using
\[ \frac{\jump{\Phi}}{\jump{h}} = \frac{g^2 \sin^2 \theta}{450 \nu^2}
   \frac{\jump{h^5}}{\jump{h}} + \frac{g \cos (\theta)}{2}
   \frac{\jump{h^2}}{\jump{h}}, \qquad \frac{\jump{\Phi'}}{\jump{h}} =
   \frac{g^2 \sin^2 \theta}{90 \nu^2} \frac{\jump{h^4}}{\jump{h}} + g \cos
   (\theta) \]
and
\[ \frac{\jump{h^2}}{\jump{h}} = h_l + h_r, \qquad \frac{\jump{h^4}}{\jump{h}}
   = (h_l + h_r) (h_l^2 + h_r^2) \]
\[ \frac{\jump{h^5}}{\jump{h}} = h_l^4 + h_l^3 h_r + h_l^2 h_r^2 + h_l h_r^3 +
   h_r^4 = (h_l + h_r) (h_l^3 + h_r^3) + h_l^2 h_r^2 \]
The computation of time average $\bar{h}$ as in \eqref{eq:barh} requires $\jump{\Phi} / \jump{h}$, which is also computed by the above formulae, but now $h_l, h_r$ represent values at two different time levels.

\section{Energy-conserving scheme for the augmented 2-D model}\label{sec:scheme2d}
The two-dimensional model \eqref{eq:aug2d}, \eqref{eq:M23} can be written in Cartesian components as
\[
\partial_t \U + \partial_1 \F_1 (\U) + \partial_2 \F_2 (\U) = \M (\U)
\]
where
\begin{equation}
  \U = \left[\begin{array}{c}
    h\\
    h v_1\\
    h v_2\\
    h w_1\\
    h w_2
  \end{array}\right], \qquad \F_1 = \left[\begin{array}{c}
    h v_1\\
    h v_1^2 + \pi (h)\\
    h v_1 v_2\\
    h v_1 w_1\\
    h v_1 w_2
  \end{array}\right], \qquad \F_2 = \left[\begin{array}{c}
    h v_2\\
    h v_1 v_2\\
    h v_2^2 + \pi (h)\\
    h v_2 w_1\\
    h v_2 w_2
  \end{array}\right], \qquad \M = \left[\begin{array}{c}
    0\\
    \bm{\M}_v\\
    \bm{\M}_w
  \end{array}\right] \label{eq:swa2d}
\end{equation}
and
\[ \bm{\M}_v
   = \left[\begin{array}{c}
     \partial_1 (h \partial_1 m_1) + \partial_2 (h \partial_1 m_2)\\
     \partial_1 (h \partial_2 m_1) + \partial_2 (h \partial_2 m_2)
   \end{array}\right] - \left[\begin{array}{c}
     \partial_1 (\gV \cdot \bm{w})\\
     \partial_2 (\gV \cdot \bm{w})
   \end{array}\right], \qquad \bm{m} = \TensV\bm{w} \]
\[ \bm{\M}_w
    = - \TensV\left[\begin{array}{c}
     \partial_1 (h \partial_1 v_1) + \partial_2 (h \partial_1 v_2)\\
     \partial_1 (h \partial_2 v_1) + \partial_2 (h \partial_2 v_2)
   \end{array}\right] - (\partial_1 v_1 + \partial_2 v_2)
   \left[\begin{array}{c}
     \gs_1\\
     \gs_2
   \end{array}\right] \]
In the above equations, $\partial_1, \partial_2$ denote partial differentiation along the $x, y$ directions, respectively, and
\begin{eqnarray*}
  q & = & \E^{- 1}_{\tmop{cap}} \left( \frac{h |\bm{w}|^2}{2 \kappa (h)}
  \right)\\
  \TensV (h, \bm{w}) & = & \sqrt{h \kappa (h)} \left[ 2 \frac{\alpha' (q^2) h}{\alpha^2
  (q^2) \kappa (h)} \bm{w} \otimes \bm{w} + \alpha (q^2) I_{2 \times 2} \right]\\
  \gV (h, \bm{w}) & = & \left[ \frac{\kappa' (h) h}{2 \kappa (h)} + \frac{1}{2} + 2
  \frac{\alpha' (q^2)}{\alpha (q^2)} q^2 \right] h \bm{w}
\end{eqnarray*}
Taking the dot product of \eqref{eq:swa2d} with entropy variables
\[ \V = \Ee' (\U)^{\top} = \left[\begin{array}{c}
     \Phi' (h) - \frac{1}{2} (| \bm{v} |^2 + | \bm{w} |^2)\\
     \bm{v}\\
     \bm{w}
   \end{array}\right] \]
we get
\[ \partial_t \Ee (\U) + \partial_1 \Gs_1 (\U) + \partial_2 \Gs_2 (\U) = \V \cdot \M \]
where
\[ \Gs_1 = (\Ee + \pi) v_1, \qquad \Gs_2 = (\Ee + \pi) v_2 \]
Integrating over the domain and assuming periodic boundary conditions
\begin{eqnarray*}
  (\V, \M) & = & \int_\Omega \V \cdot \M \mathd x\\
  &=& (\tmop{div} (h \nabla (\TensV\bm{w})^{\top}), \bm{v}) -
  (\nabla (\gV \cdot \bm{w}), \bm{v}) - (\TensV\tmop{div} (h
  \nabla \bm{v}^{\top}), \bm{w}) - (\gV \tmop{div}
  \bm{v}, \bm{w})\\
  & = & (\tmop{div} (h \nabla (\TensV\bm{w})^{\top}), \bm{v}) -
  (\tmop{div} (h \nabla \bm{v}^{\top}), \TensV\bm{w}) - (\nabla
  (\gV \cdot \bm{w}), \bm{v}) - (\tmop{div} \bm{v},
  \gV \cdot \bm{w})\\
  & = & \int_{\Omega} \tmop{div} (\nabla (\TensV\bm{w}) h \bm{v} - h
  (\nabla \bm{v}) \TensV\bm{w} - (\gV \cdot \bm{w})
  \bm{v}) \mathd x\\
  & = & \int_{\partial \Omega} (\nabla (\TensV\bm{w}) h \bm{v} - h
  (\nabla \bm{v}) \TensV\bm{w} - (\gV \cdot \bm{w})
  \bm{v}) \cdot \bm{n} \mathd s \\
  &=& 0
\end{eqnarray*}
Consider a partition of the rectangular domain into cells of size $(\Delta x,
\Delta y)$ indexed by $(j, k)$. The fully discrete scheme is given by
\begin{equation}
  \frac{\U_{j, k}^{n + 1} - \U_{j, k}^n}{\Delta t} + \frac{(\F_1^{\star})_{\jph, k}^{\nph} - (\F_1^{\star})_{\jmh, k}^{\nph}}{\Delta x} +
  \frac{(\F_2^{\star})_{j, \kph}^{\nph} - (\F_2^{\star})_{j, \kmh}^{\nph}}{\Delta y} = \M_{j, k}^{\nph}\label{eq:fd2d}
\end{equation}
where $(\F_1^{\star})_{\jph, k} = \F_1^{\star} (\U_{j, k}, \U_{j + 1, k})$,  $(\F_2^{\star})_{j, \kph} = \F_2^{\star} (\U_{j, k}, \U_{j, k + 1})$ and $\M_{j, k}^{\nph} = \M_{j, k} (\U^{\nph})$. For
the case of $\Phi (h) = \frac{1}{2} g h^2$, these fluxes are given by
\[ \F_1^{\star} = \left[\begin{array}{c}
     \avgc{h} \avgc{v_1}\\[2mm]
     \avgc{h} \avgc{v_1}^2 + \frac{1}{2} g \avgc{h^2}\\[2mm]
     \avgc{h} \avgc{v_1} \avgc{v_2}\\[2mm]
     \avgc{h} \avgc{v_1} \avgc{w_1}\\[2mm]
     \avgc{h} \avgc{v_1} \avgc{w_2}
   \end{array}\right], \qquad \F_2^{\star} = \left[\begin{array}{c}
     \avgc{h} \avgc{v_2}\\[2mm]
     \avgc{h} \avgc{v_1} \avgc{v_2}\\[2mm]
     \avgc{h} \avgc{v_2}^2 + \frac{1}{2} g \avgc{h^2}\\[2mm]
     \avgc{h} \avgc{v_2} \avgc{w_1}\\[2mm]
     \avgc{h} \avgc{v_2} \avgc{w_2}
   \end{array}\right] \]
and for general $\Phi(h)$ by
\[
  \F^{\star}_1 =
    \left[\begin{array}{c}
    \Fs_{11}^{\star} \\[2mm]
   \Fs_{21}^{\star} \\[2mm]
   \Fs_{31}^{\star} \\[2mm]
    \Fs_{41}^{\star} \\[2mm]
    \Fs_{51}^{\star}
  \end{array}\right]  =
  \left[\begin{array}{c}
    \frac{\jump{\pi}}{\jump{\Phi'}} \avgc{v_1}\\[2mm]
    \avgc{v_1} \Fs_{11}^{\star} + \avgc{\pi}\\[2mm]
    \avgc{v_2} \Fs_{11}^{\star}\\[2mm]
    \avgc{w_1} \Fs_{11}^{\star}\\[2mm]
    \avgc{w_2} \Fs_{11}^{\star}
  \end{array}\right], \qquad
  \F^{\star}_2 =
    \left[\begin{array}{c}
    \Fs_{12}^{\star} \\[2mm]
   \Fs_{22}^{\star} \\[2mm]
   \Fs_{32}^{\star} \\[2mm]
    \Fs_{42}^{\star} \\[2mm]
    \Fs_{52}^{\star}
  \end{array}\right]  =
  \left[\begin{array}{c}
    \frac{\jump{\pi}}{\jump{\Phi'}} \avgc{v_2}\\[2mm]
    \avgc{v_1} \Fs_{12}^{\star}\\[2mm]
    \avgc{v_2} \Fs_{12}^{\star} + \avgc{\pi}\\[2mm]
    \avgc{w_1} \Fs_{12}^{\star}\\[2mm]
    \avgc{w_2} \Fs_{12}^{\star}
  \end{array}\right]
\]
The time-averaged state $\U^{\nph}$ is designed to satisfy \eqref{eq:tavg}.  In terms of primitive variables, the time average is given by
\[
\bar{\bm{v}} = \avgc{\bm{v}}, \qquad
\bar{\bm{w}} = \avgc{\bm{w}}
\]
and $\bar h$ is obtained from
\[
\Phi'(\bar h) = \frac{\jump{\Phi}}{\jump{h}} + \frac{1}{8} |\jump{\bm{v}}|^2 + \frac{1}{8} |\jump{\bm{w}}|^2
\]
In the case of $\Phi(h) = \frac{1}{2} g h^2$, we can write $\bar h$ explicitly as
\[
\bar h = \avgc{h} + \frac{1}{8g} |\jump{\bm{v}}|^2 + \frac{1}{8g} |\jump{\bm{w}}|^2
\]
For the discretization of the capillary terms $\M$, define the finite difference and averaging operations along the coordinate directions by
\[ \delta_1 a_{j, k} = \frac{a_{\jph, k} - a_{\jmh, k}}{\Delta x},
   \qquad \mu_1 a_{\jph, k} = \frac{1}{2} (a_{j, k} + a_{j + 1, k}) \]
\[ \delta_2 a_{j, k} = \frac{a_{j, \kph} - a_{j, \kmh}}{\Delta y},
   \qquad \mu_2 a_{j, \kph} = \frac{1}{2} (a_{j, k} + a_{j, k + 1}) \]
Then
\[ (\bm{\M}_v)_{j, k} = \left[\begin{array}{c}
     \delta_1 [(\mu_1 h) \delta_1 m_1] + \delta_2 \mu_2 (h \delta_1 \mu_1
     m_2)\\
     \delta_1 \mu_1 (h \delta_2 \mu_2 m_1) + \delta_2 [(\mu_2 h) \delta_2 m_2]
   \end{array}\right]_{j, k} - \left[\begin{array}{c}
     \delta_1 \mu_1 (\gV \cdot \bm{w})\\
     \delta_2 \mu_2 (\gV \cdot \bm{w})
   \end{array}\right]_{j, k} \]
and
\[ (\bm{\M}_w)_{j, k} = - \TensV_{j, k} \left[\begin{array}{c}
     \delta_1 [(\mu_1 h) \delta_1 v_1] + \delta_2 \mu_2 (h \delta_1 \mu_1
     v_2)\\
     \delta_1 \mu_1 (h \delta_2 \mu_2 v_1) + \delta_2 [(\mu_2 h) \delta_2 v_2]
   \end{array}\right]_{j, k} - (\delta_1 \mu_1 v_1 + \delta_2 \mu_2 v_2)_{j,
   k} \left[\begin{array}{c}
     \gs_1\\
     \gs_2
   \end{array}\right]_{j, k} \]
Following Bresch et al.~\cite{Bresch2020}, we can show that
\[ 
(\V, \M)_\Delta := \Delta x \Delta y \sum_{j, k} \V_{j, k} \cdot \M_{j, k} = \Delta x
   \Delta y \sum_{j, k} [\bm{v}_{j, k} \cdot (\bm{\M}_v)_{j, k} +
   \bm{w}_{j, k} \cdot (\bm{\M}_w)_{j, k}] = 0 
\]
under periodic boundary conditions. In fact
\begin{eqnarray}
  &  & (\bm{\M}_v, \bm{v})_\Delta + (\bm{\M}_w, \bm{w})_\Delta
  \nonumber\\
  & = & \left( \left[\begin{array}{c}
    \delta_1 [(\mu_1 h) \delta_1 m_1] + \delta_2 \mu_2 (h \delta_1 \mu_1
    m_2)\\
    \delta_1 \mu_1 (h \delta_2 \mu_2 m_1) + \delta_2 [(\mu_2 h) \delta_2 m_2]
  \end{array}\right], \bm{v} \right)_\Delta - \left( \left[\begin{array}{c}
    \delta_1 \mu_1 (\gV \cdot \bm{w})\\
    \delta_2 \mu_2 (\gV \cdot \bm{w})
  \end{array}\right], \bm{v} \right)_\Delta \nonumber\\
  &  & - \left( \TensV\left[\begin{array}{c}
    \delta_1 [(\mu_1 h) \delta_1 v_1] + \delta_2 \mu_2 (h \delta_1 \mu_1
    v_2)\\
    \delta_1 \mu_1 (h \delta_2 \mu_2 v_1) + \delta_2 [(\mu_2 h) \delta_2 v_2]
  \end{array}\right], \bm{w} \right)_\Delta - \left( (\delta_1 \mu_1 v_1 +
  \delta_2 \mu_2 v_2)  \left[\begin{array}{c}
    \gs_1\\
    \gs_2
  \end{array}\right], \bm{w} \right)_\Delta \nonumber\\
  & = & \left( \left[\begin{array}{c}
    \delta_1 [(\mu_1 h) \delta_1 m_1] + \delta_2 \mu_2 (h \delta_1 \mu_1
    m_2)\\
    \delta_1 \mu_1 (h \delta_2 \mu_2 m_1) + \delta_2 [(\mu_2 h) \delta_2 m_2]
  \end{array}\right], \bm{v} \right)_\Delta - \left( \left[\begin{array}{c}
    \delta_1 [(\mu_1 h) \delta_1 v_1] + \delta_2 \mu_2 (h \delta_1 \mu_1
    v_2)\\
    \delta_1 \mu_1 (h \delta_2 \mu_2 v_1) + \delta_2 [(\mu_2 h) \delta_2 v_2]
  \end{array}\right], \bm{m} \right)_\Delta \nonumber\\
  &  & - \left( \left[\begin{array}{c}
    \delta_1 \mu_1 (\gV \cdot \bm{w})\\
    \delta_2 \mu_2 (\gV \cdot \bm{w})
  \end{array}\right], \bm{v} \right)_\Delta - ((\delta_1 \mu_1 v_1 + \delta_2
  \mu_2 v_2), \gV \cdot \bm{w})_\Delta  \label{eq:Msbp}
\end{eqnarray}
The first two terms on the right of \eqref{eq:Msbp} are equal by Lemma 3.2 of~\cite{Bresch2020}; in fact this follows from \eqref{eq:sbp2} and an additional SBP relation~\eqref{eq:sbp2da}. The last two terms of \eqref{eq:Msbp} are
\[ 
- [(\delta_1 \mu_1 (\gV \cdot \bm{w}), v_1)_\Delta + (\delta_1
   \mu_1 v_1, \gV \cdot \bm{w})_\Delta] - [(\delta_2 \mu_2
   (\gV \cdot \bm{w}), v_2)_\Delta + (\delta_2 \mu_2 v_2,
   \gV \cdot \bm{w})_\Delta] = 0 
\]
which follows from \eqref{eq:sbp1}. Taking the dot product of the fully discrete
scheme \eqref{eq:fd2d} with $\V_{j, k}^{\nph} = \V (\bar{\U} (\U_{j, k}^n,
\U_{j, k}^{n + 1}))$, we get
\begin{equation}
     \frac{\Ee (\U_{j, k}^{n + 1}) - \Ee (\U_{j, k}^n)}{\Delta t}  + \displaystyle
     \frac{(\Gs_1^{\star})_{\jph, k}^{\nph} - (\Gs_1^{\star})_{\jmh, k}^{\nph}}{\Delta x}  + \displaystyle\frac{(\Gs_2^{\star})_{j, \kph}^{\nph} - (\Gs_2^{\star})_{j, \kmh}^{\nph}}{\Delta y} =  \displaystyle\V_{j, k}^{\nph} \cdot \M_{j, k}^{\nph}
\end{equation}
where
\[
\Gs_m^*(\U_l, \U_r) = \left( \frac{\V_l + \V_r}{2} \right) \cdot
   \F^{\star}_m (\U_l, \U_r) - \frac{\H_m(\U_l) + \H_m(\U_r)}{2}, \quad m=1,2
\]
where $\Gs^*_m$ is a consistent numerical energy flux and $\H_m = \V \cdot \F_m - \Gs_m$. By summation over all cells, we obtain that the total energy is conserved under periodic boundary conditions.

\subsection{Application to thin-film flow model}
The model for two-dimensional thin-film flow down an inclined plane at an angle $\theta$, which is an extension of the one-dimensional model from Section~\ref{sec:richard}, can be written as
\[
\partial_t \U + \partial_1 \F_1 (\U) + \partial_2 \F_2 (\U) = \M (\U) + \D (\U)
\]
where the various terms are given in the last section, and $\D (\U)$ contains the frictional terms and the gravitational force, and is given by
\[
\D (\U) = \left[\begin{array}{c}
     0\\[1mm]
     g h \sin \theta - \frac{3 \nu v_1}{h} + \frac{9 \nu}{2} \partial_1 (h
     \tmop{div} \bm{v})\\[2mm]
     - \frac{3 \nu v_2}{h} + \frac{9 \nu}{2} \partial_2 (h \tmop{div}
     \bm{v})\\[1mm]
     0\\[1mm]
     0
   \end{array}\right]
\]
The pressure and potential functions are given by \eqref{eq:richpi}, \eqref{eq:richPhi}. The parabolic terms are approximated in a multi-dimensional manner,
\begin{equation*}
  \partial_1 (h \tmop{div} \bm{v})_{j, k}  \approx  \delta_1 \mu_2 (h
  \nobracket \tmop{div} \bm{v})_{j, k}, \qquad
  \partial_2 (h \tmop{div} \bm{v})_{j, k}  \approx  \delta_2 \mu_1 (h
  \nobracket \tmop{div} \bm{v})_{j, k}
\end{equation*}
which requires quantities at vertices that are defined as
\begin{equation}\label{eq:hdivvert}
\begin{aligned}
  h_{\jph, \kph} & =  (\mu_1 \mu_2 h)_{\jph, \kph}\\
  (\tmop{div} \bm{v})_{\jph, \kph} & =  (\delta_1 \mu_2 v_1 +
  \delta_2 \mu_1 v_2)_{\jph, \kph}
\end{aligned}
\end{equation}
With this approximation, the contribution of the parabolic terms to the total energy equation is
\begin{eqnarray*}
\E_D &=& \frac{9\nu}{2} [ (\delta_1 \mu_2 (h \nobracket \tmop{div} \bm{v}) \nobracket, v_1)_\Delta +
   (\delta_2 \mu_1 (h \nobracket \tmop{div} \bm{v}) \nobracket, v_2)_\Delta] \\
   &=& - \frac{9\nu}{2} \sum_{j, k} h_{\jph, \kph} (\tmop{div}
   \bm{v})_{\jph, \kph}^2 \Delta x \Delta y  \\
   &\le& 0
\end{eqnarray*}
where we used the SBP relation~\eqref{eq:sbpdiv}  assuming periodic boundary conditions and positivity of $h$. The frictional terms also make a non-positive contribution to the right hand side of the total energy equation so that they do not increase the total energy.

\begin{remark}
We can approximate the parabolic terms in a simpler one-dimensional manner as
\begin{eqnarray*}
  \partial_1 (h \tmop{div} \bm{v})_{j, k} & \approx & \delta_1 ((\mu_1
  h)  (\tmop{div} \bm{v}))_{j, k}\\
  \partial_2 (h \tmop{div} \bm{v})_{j, k} & \approx & \delta_2 ((\mu_2
  h)  (\tmop{div} \bm{v}))_{j, k}
\end{eqnarray*}
which requires $\tmop{div} \bm{v}$ at the faces computed as
\begin{eqnarray*}
  (\tmop{div} \bm{v})_{\jph, k} & \approx & (\delta_1 v_1 +
  \delta_2 \mu_1 \mu_2 v_2)_{\jph, k}\\
  (\tmop{div} \bm{v})_{j, \kph} & \approx & (\delta_1 \mu_1 \mu_2
  v_1 + \delta_2 v_2)_{j, \kph}
\end{eqnarray*}
But in this case, we are not able to show the dissipative character of these
terms. We obtain
\begin{eqnarray*}
  &  & (\delta_1 ((\mu_1 h)  (\tmop{div} \bm{v})), v_1)_\Delta + (\delta_2
  ((\mu_2 h)  (\tmop{div} \bm{v})), v_2)_\Delta\\
  & = & - \Delta x \Delta y \sum_{j, k} [ (\mu_1 h) (\delta_1 v_1)^2]_{\jph, k} - \Delta x \Delta y \sum_{j, k} [(\mu_2 h) (\delta_2 v_2)^2]_{j, \kph}  \\
  &  & + \Delta x \Delta y \sum_{j, k} (v_1)_{j, k} \delta_1 ((\mu_1 h)
  (\delta_2 \mu_1 \mu_2 v_2))_{j, k}\\
  &  & + \Delta x \Delta y \sum_{j, k} (v_2)_{j, k} \delta_2 ((\mu_2 h)
  (\delta_1 \mu_1 \mu_2 v_1))_{j, k}
\end{eqnarray*}
but the sign of the last two terms is indeterminate. Thus a good choice of the finite difference operators is necessary to obtain the correct sign of the parabolic terms in the energy equation.
\end{remark}

\section{Numerical results}\label{sec:numerics}
The proposed method is implemented in a C code that uses the Message Passing Interface (MPI) for parallelization and PETSc for solving linear and nonlinear equations. We solve the nonlinear implicit system generated at each time step by the Crank-Nicolson scheme using Newton's method coupled with a matrix-free GMRES Krylov subspace solver. We use no preconditioner for GMRES since the time steps are already small. The time step is chosen based on a hyperbolic CFL condition, which in one dimension is of the form
\[ \Delta t = \tmop{CFL} \cdot \frac{\Delta x}{\max_j \left( | v_j | +
   \sqrt{\pi' (h_j)} \right)} \]
Unless stated otherwise, we use $\tmop{CFL}=0.1$ for all test cases. All computations are performed in dimensional units.

We will assume that we know the spatial derivatives of the initial height field. Then, under the assumption that $\kappa (h) = \sigma/ \rho$ with $\sigma$ constant, the initial condition for $w$ is set as
\[
w = w (h, \partial_x h) = \sqrt{\frac{2 \sigma\E_{\tmop{cap}}
   (| \partial_x h |)}{\rho h}} \tmop{sign} (\partial_x h)
\]
since $w$ must have the same sign as $\partial_x h$. We use a similar strategy in the two-dimensional case, leading to
\begin{equation}
  \bm{w} = \sqrt{\frac{2 \sigma\E_{\tmop{cap}} (| \nabla h
  |)}{\rho h}}  \frac{\nabla h}{| \nabla h |} \label{eq:wgradh}
\end{equation}
For constant $\kappa(h)$, we list other expressions required for the implementation of this model in a computer code.

\paragraph{Quadratic capillary model.} The capillary energy function is
\[ \E_{\tmop{cap}} (q) = \frac{1}{2} q^2, \qquad \alpha (s) =
   \sqrt{\frac{2\E_{\tmop{cap}} \left( \sqrt{s} \right)}{s}} = 1 \]
Then in one dimension, we have
\[ \TensS(h, w) = \sqrt{\frac{h \sigma}{\rho}}, \qquad \gs (h, w) = \frac{1}{2} h w
\]
and in two dimensions, they are given by
\[ \TensV(h, \bm{w}) = \sqrt{\frac{h \sigma}{\rho}} I_{2 \times 2}, \qquad \gV
   (h, \bm{w}) = \frac{1}{2} h \bm{w} \]
\paragraph{Nonlinear capillary model.} The capillary energy function is
\[ \E_{\tmop{cap}} (q) = \sqrt{1 + q^2} - 1, \qquad \alpha (s) =
   \sqrt{\frac{2}{1 + \sqrt{1 + s}}} \]
Then in one dimension, we have
\[
  \TensS(h, w)  =  \sqrt{\frac{\sigma h / \rho}{1 + \frac{\rho h w^2}{4
  \sigma}}} \left( 1 - \frac{\rho h w^2}{4 \sigma \left( 1 + \frac{\rho h
  w^2}{2 \sigma} \right)} \right), \qquad
  \gs (h, w) =  \frac{h w}{2 \left( 1 + \frac{\rho h w^2}{2 \sigma} \right)}
\]
and in two dimensions, we have
\[ \TensV(h, \bm{w}) = \sqrt{\frac{\sigma h / \rho}{1 + \frac{\rho h |
   \bm{w} |^2}{4 \sigma}}} \left( I_{2 \times 2} - \frac{\rho h
   \bm{w} \otimes \bm{w}}{4 \sigma \left( 1 + \frac{\rho h |
   \bm{w} |^2}{2 \sigma} \right)} \right), \qquad \gV (h, \bm{w}) =
   \frac{h \bm{w}}{2 \left( 1 + \frac{\rho h | \bm{w} |^2}{2
   \sigma} \right)} \]

\subsection{Gaussian bump in 1-D}\label{sec:gauss1d}

\begin{figure}[htbp]
    \centering
    \begin{subfigure}[t]{0.48\textwidth}
        \centering
        \includegraphics[width=\textwidth]{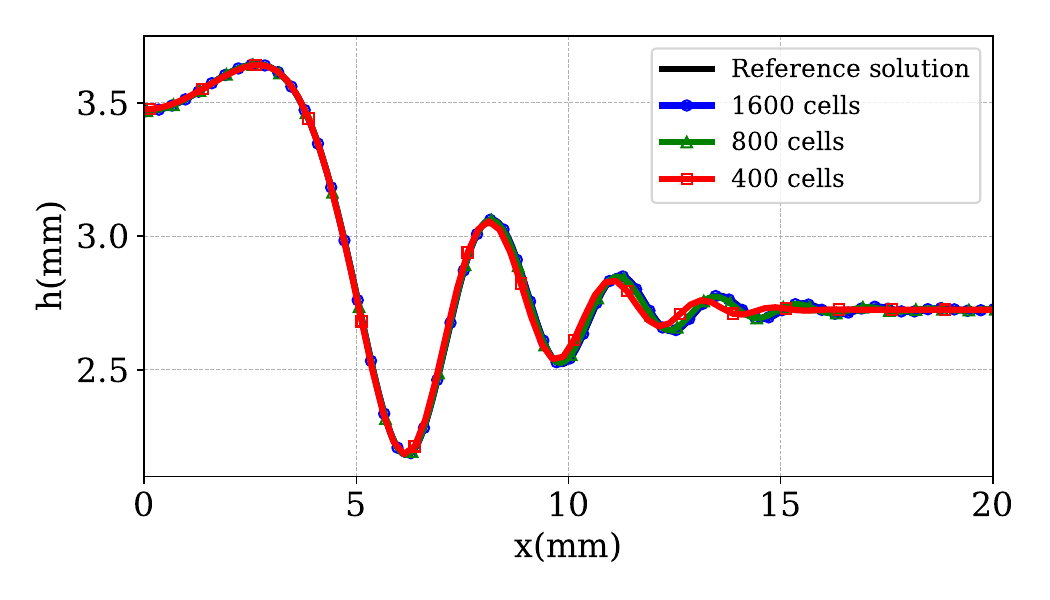}
    \end{subfigure}
    \hfill
    \begin{subfigure}[t]{0.48\textwidth}
        \centering
        \includegraphics[width=\textwidth]{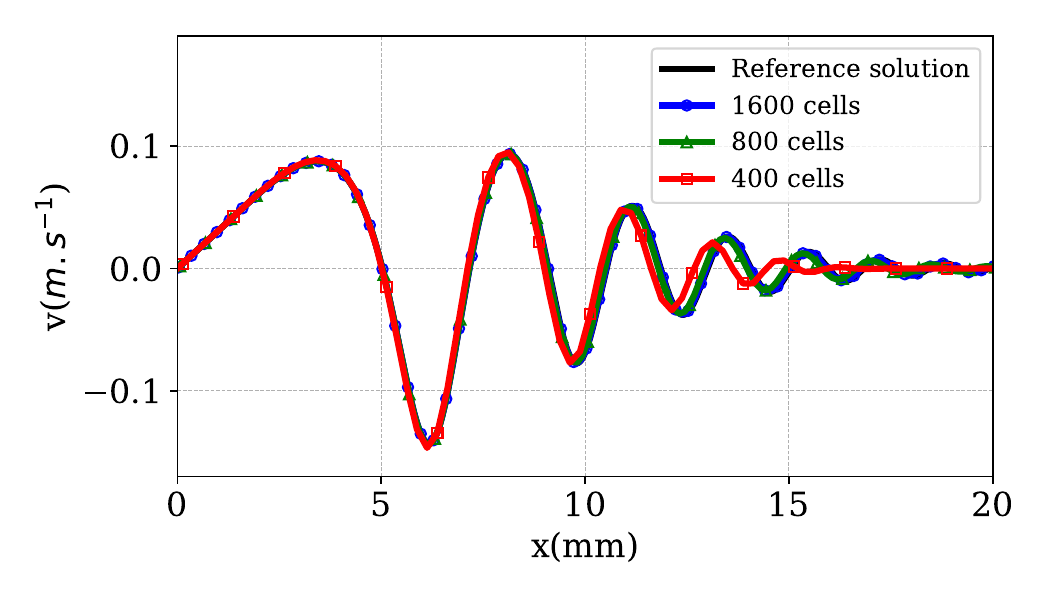}
    \end{subfigure}
    \caption{1-D Gaussian bump: grid convergence for the quadratic capillary model. The reference solution is computed with 10,000 cells.}
    \label{fig:g1dquadnx}
\end{figure}

\begin{figure}[htbp]
    \centering
    \begin{subfigure}[t]{0.48\textwidth}
        \centering
        \includegraphics[width=\textwidth]{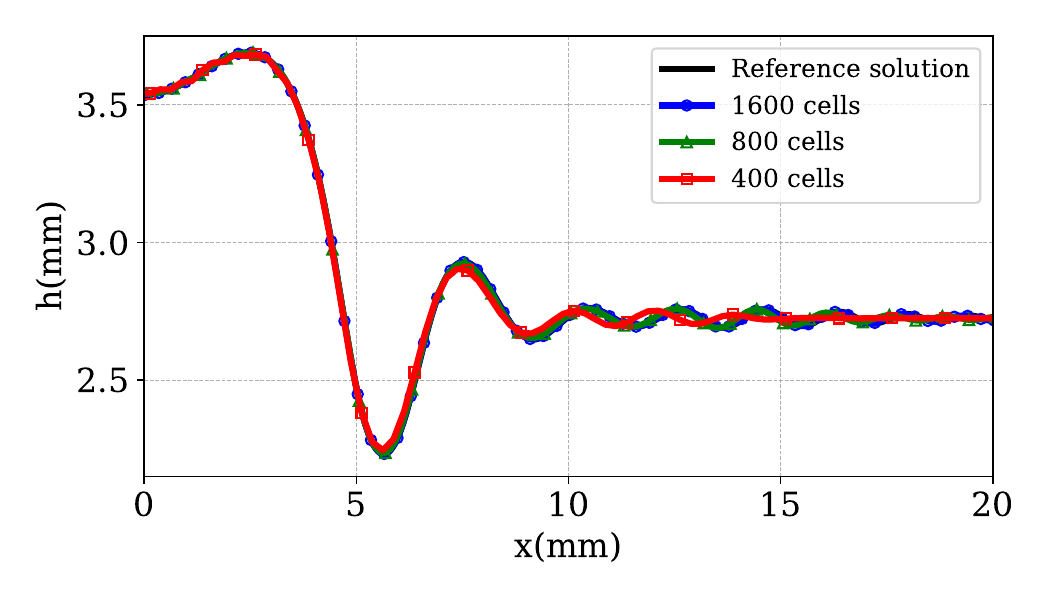}
    \end{subfigure}
    \hfill
    \begin{subfigure}[t]{0.48\textwidth}
        \centering
        \includegraphics[width=\textwidth]{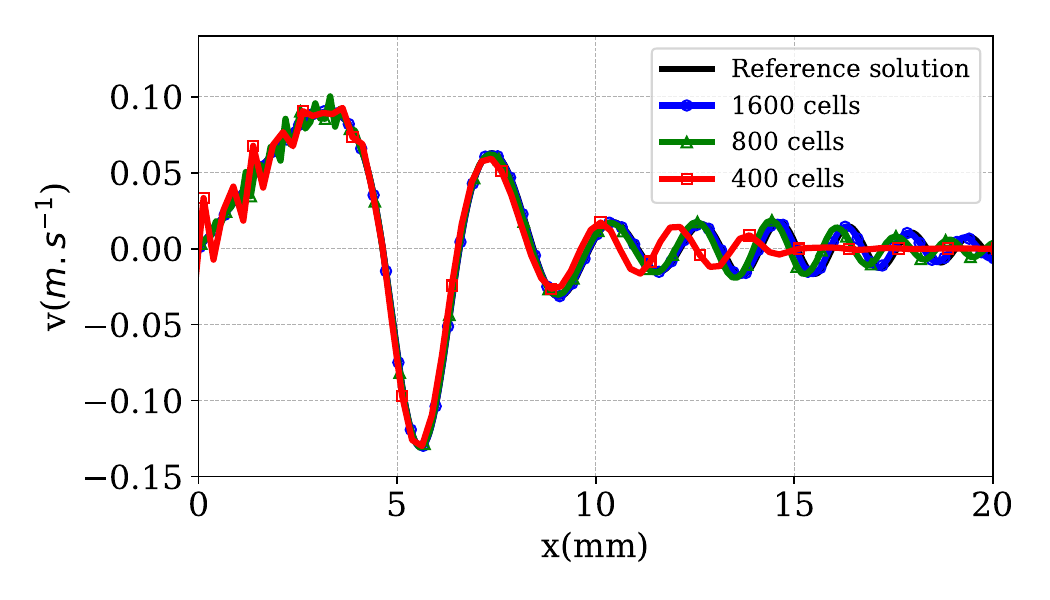}
    \end{subfigure}
    \caption{1-D Gaussian bump: grid convergence for the nonlinear capillary model. The reference solution is computed with 10,000 cells.}
    \label{fig:g1dnonnx}
\end{figure}

\begin{figure}[htbp]
    \centering
        \begin{subfigure}[t]{0.48\textwidth}
        \centering
        \includegraphics[width=\textwidth]{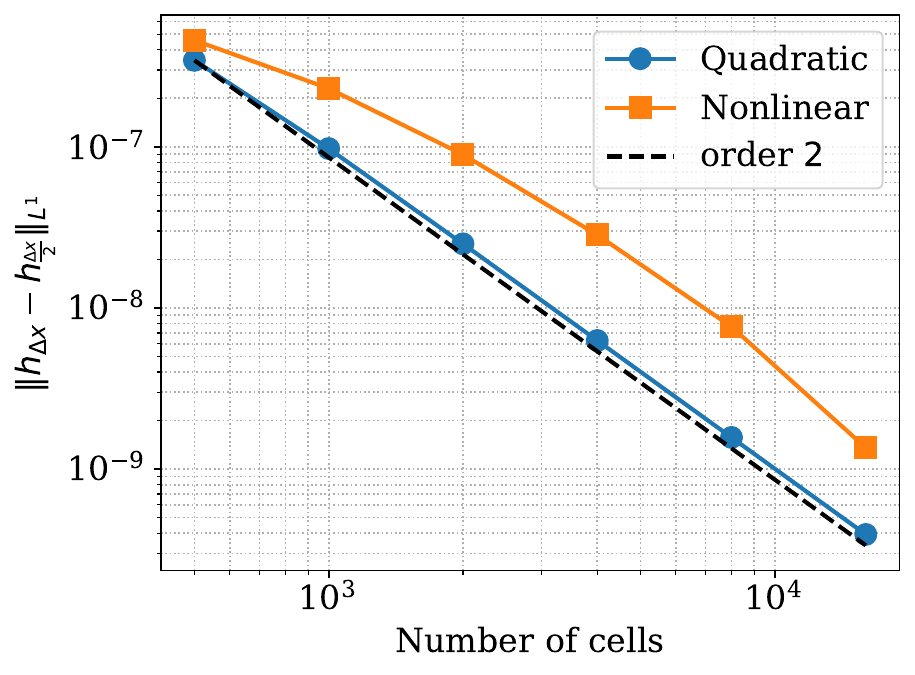}
    \end{subfigure}
    \hfill
    \begin{subfigure}[t]{0.48\textwidth}
        \centering
        \includegraphics[width=\textwidth]{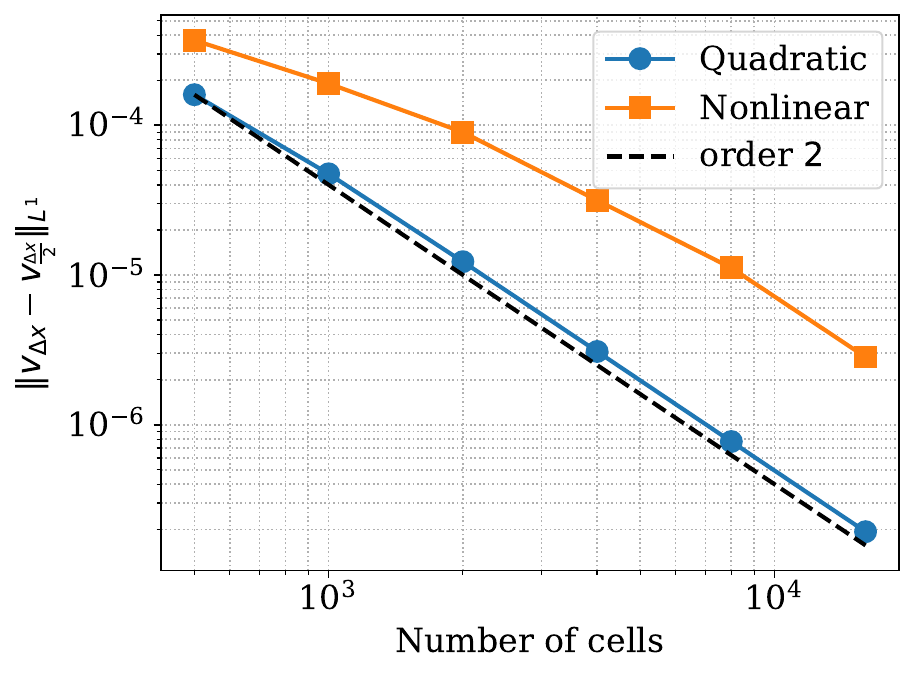}
    \end{subfigure}
    \caption{1-D Gaussian bump: convergence under grid refinement}
    \label{fig:g1dconv}
\end{figure}

This test case is taken from~\cite{Bresch2020} and consists of an initial Gaussian-shaped free surface which deforms under gravitational and capillary forces. The domain is $[- 50 \tmop{mm}, 50 \tmop{mm}]$ with periodic boundary conditions and the solution is computed up to time $5 \times 10^{-3}$ seconds. The initial condition is given by
\[ h (x, 0) = h_0 + h_1 \exp \left( - \frac{x^2}{2 (b / b_0)^2} \right),
   \qquad v (x, 0) = 0 \]
where
\[
h_1 = 2.725 \tmop{mm}, \qquad h_0 = h_1, \qquad b = 1.5 h_1, \qquad b_0 =
   4.29193
\]
The other parameters are
\[
\sigma = 0.0728~\tmop{N}/\tmop{m}, \qquad \rho = 1000~\tmop{kg}/\tmop{m}^3, \qquad g = 9.81~\tmop{m}/\tmop{s}^2
\]
Figures~\ref{fig:g1dquadnx},~\ref{fig:g1dnonnx} show the solution at the final time for the quadratic and nonlinear capillary models, respectively. The coarse grids do not resolve all the waves but the resolution improves with grid refinement. The nonlinear capillary case shows spurious numerical oscillations on coarse grids since we do not have any implicit numerical dissipation in our scheme. However, these oscillations vanish when the grid is refined. Figure~\ref{fig:g1dconv} shows the convergence of the $L^1$ error norm in $h$ under grid refinement for both models. Since we do not have an exact solution, we plot the norm of the solution difference obtained with grid sizes $\Delta x$ and $\Delta x/2$. In the case of the quadratic capillary model, a second-order rate is clearly observed on all grid levels, but this is hard to achieve for the nonlinear model because of the small-scale oscillations that require very fine grids. However, asymptotically we observe close to a second-order rate even for the nonlinear capillary model.

\begin{figure}[htb]
    \centering
    \begin{subfigure}[t]{0.48\textwidth}
        \centering
        \includegraphics[width=\textwidth]{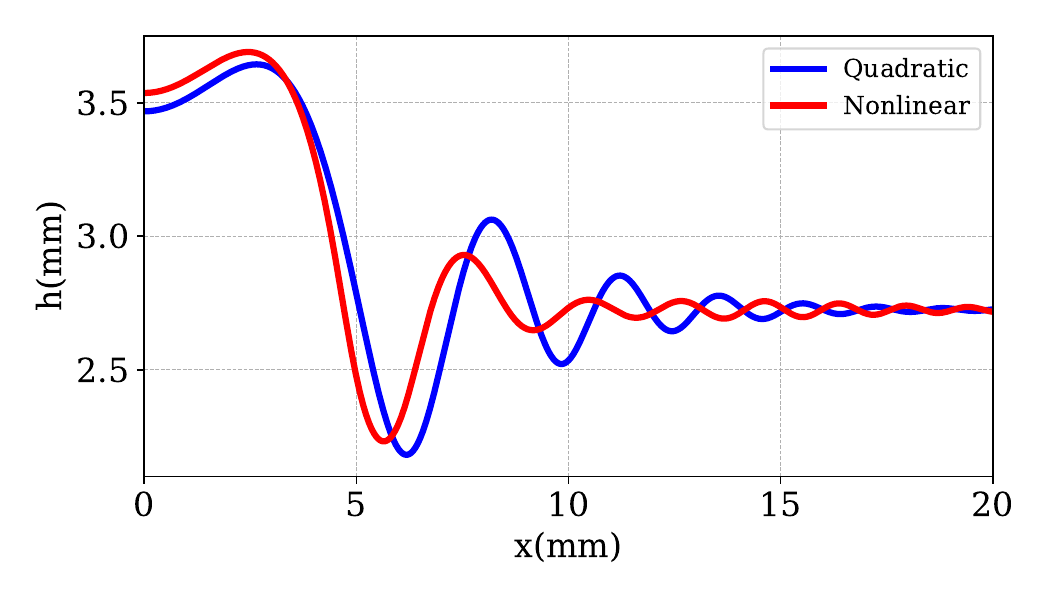}
    \end{subfigure}
    \hfill
    \begin{subfigure}[t]{0.48\textwidth}
        \centering
        \includegraphics[width=\textwidth]{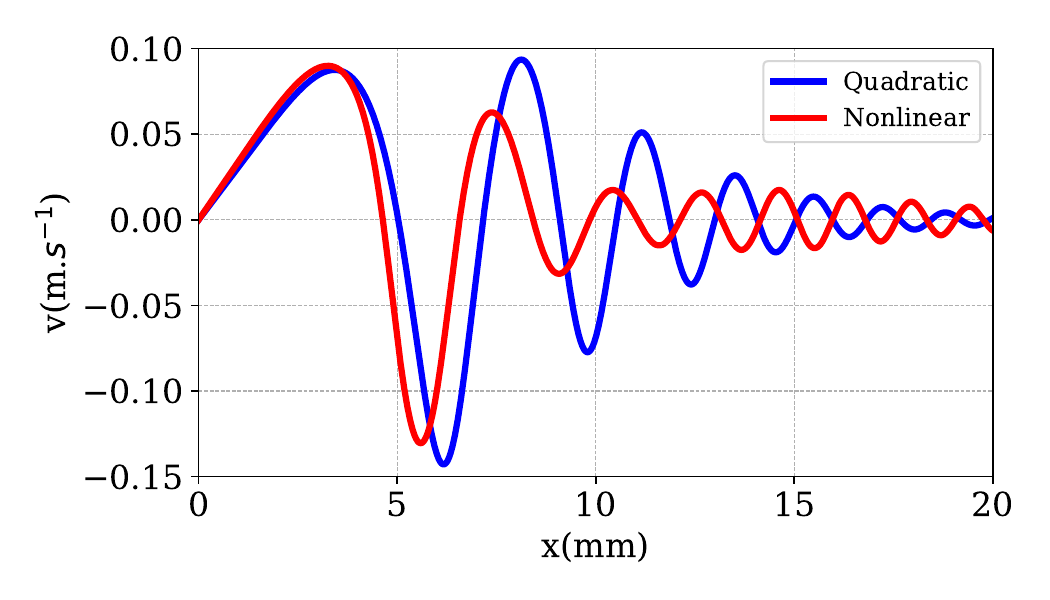}
    \end{subfigure}\\
      \begin{subfigure}[t]{0.48\textwidth}
        \centering
        \includegraphics[width=\textwidth]{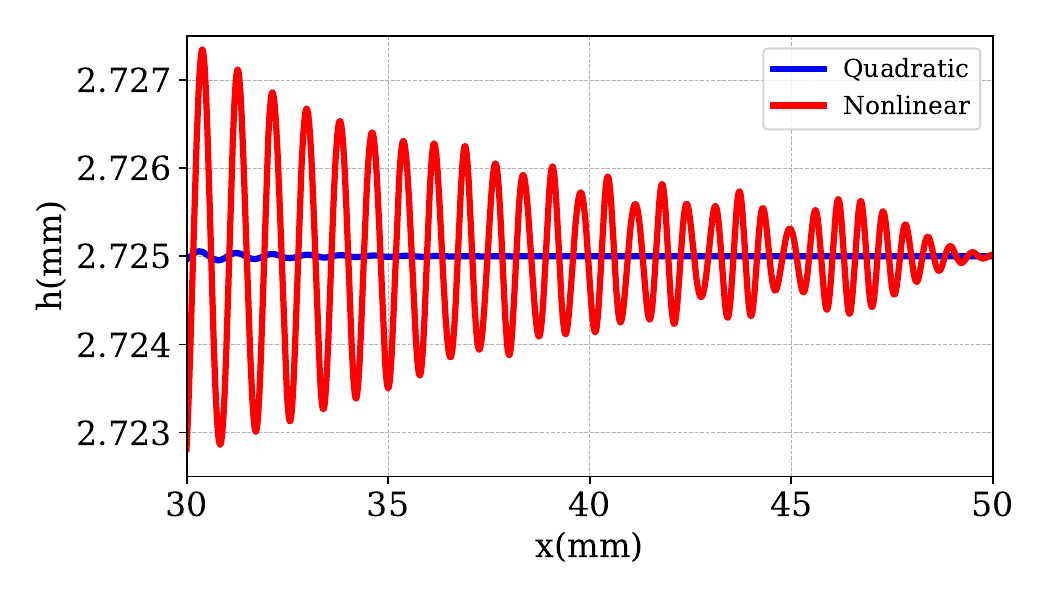}
    \end{subfigure}
    \hfill
    \begin{subfigure}[t]{0.48\textwidth}
        \centering
        \includegraphics[width=\textwidth]{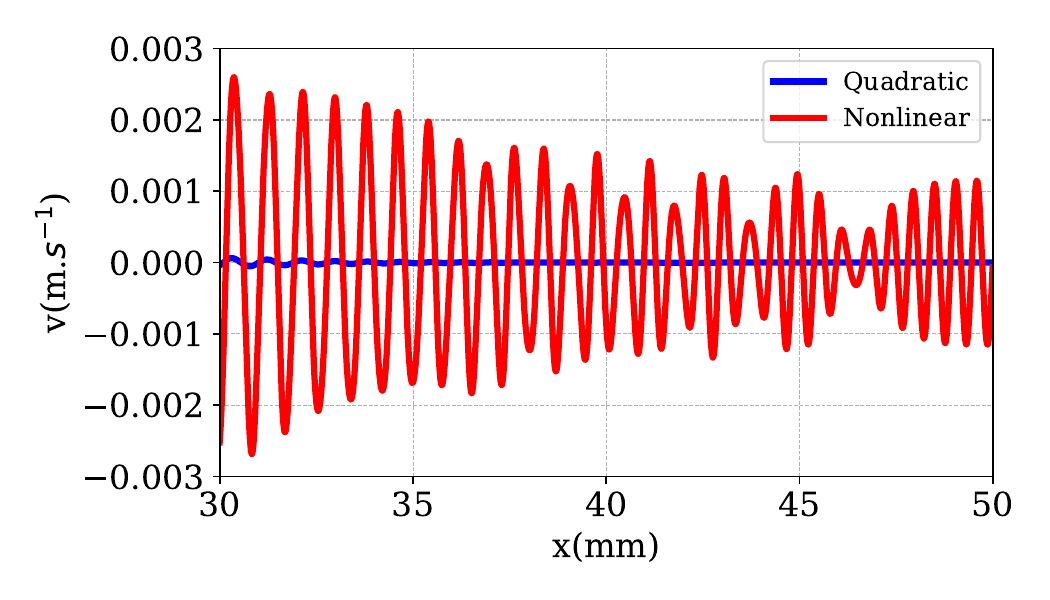}
    \end{subfigure}
    \caption{1-D Gaussian bump: comparison of the quadratic and nonlinear capillary models using 10,000 cells}
    \label{fig:g1d10k}
\end{figure}

Figure~\ref{fig:g1d10k} shows the solutions on a grid of 10,000 cells, with the second row showing the solution in the region $x\in[30 \tmop{mm}, 50 \tmop{mm}]$. The difference between the quadratic and nonlinear capillary models is clear, with the nonlinear model having attenuated waves near the center. However, there is significantly more amplitude in the nonlinear case far away from the center.

\begin{figure}[htbp]
    \centering
    \begin{subfigure}[t]{0.48\textwidth}
        \centering
        \includegraphics[width=\textwidth]{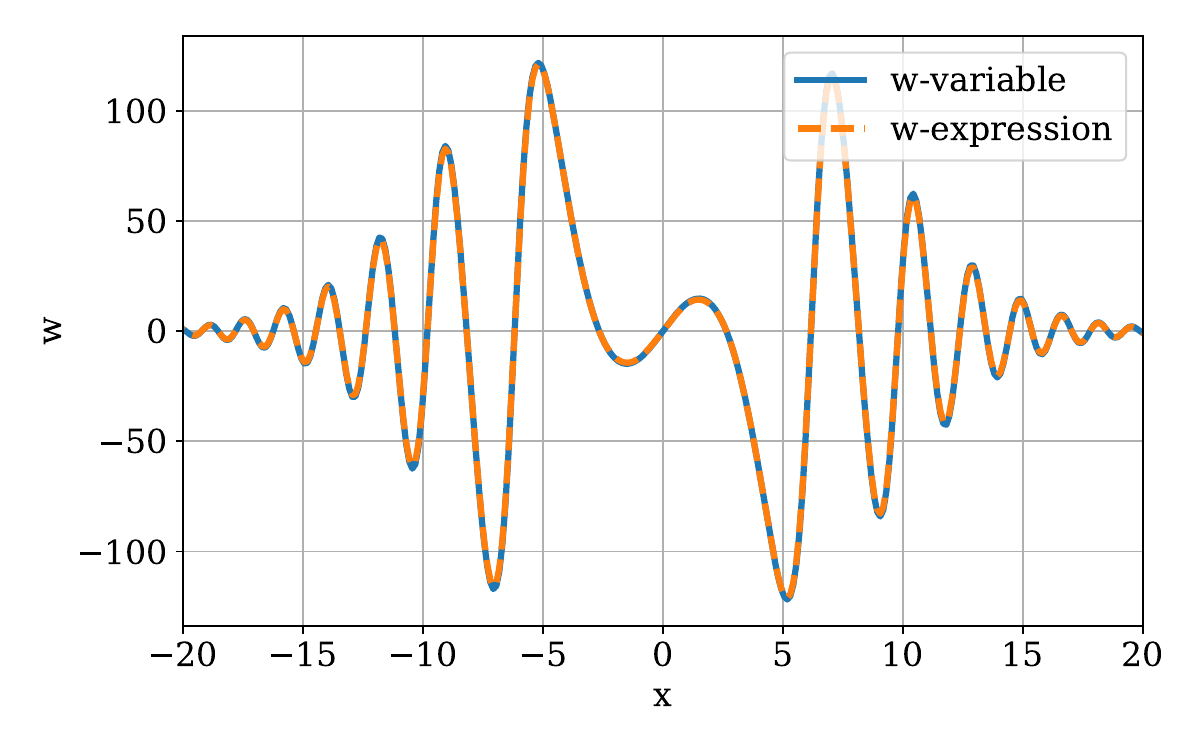}
        \caption{800 cells}
    \end{subfigure}
    \hfill
    \begin{subfigure}[t]{0.48\textwidth}
        \centering
        \includegraphics[width=\textwidth]{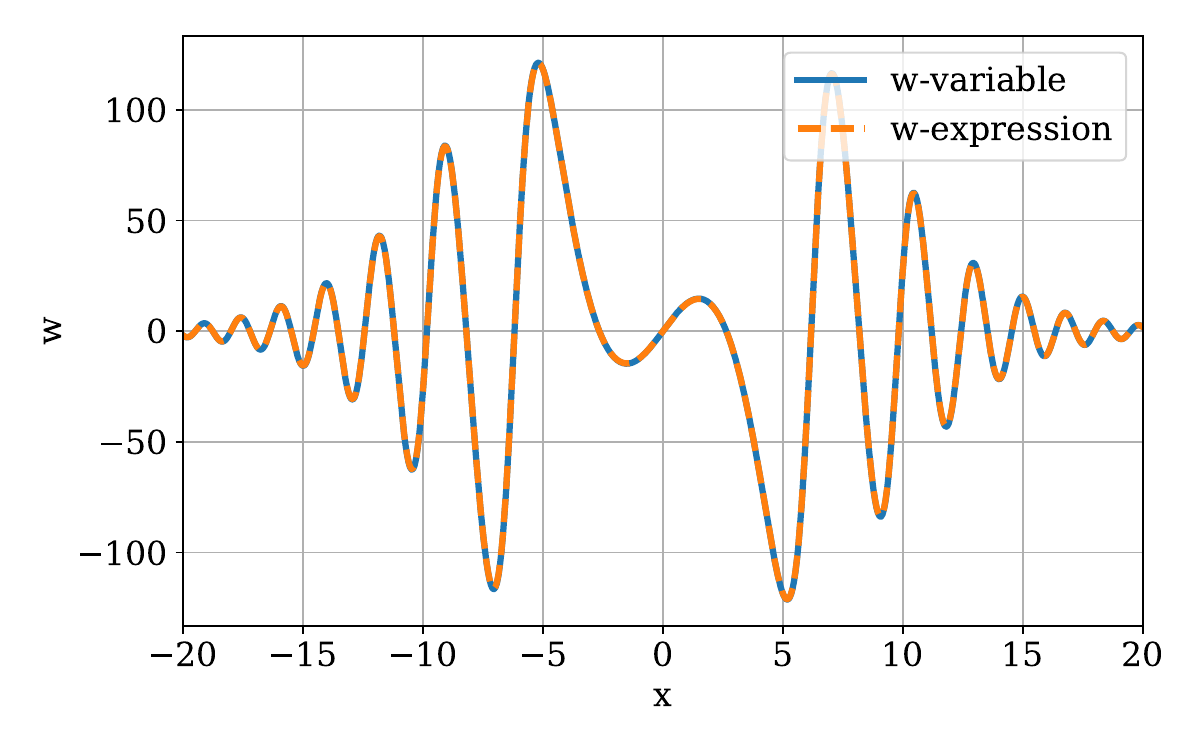}
        \caption{1600 cells}
    \end{subfigure}
    \caption{1-D Gaussian bump: comparison of  $w$ vs $\sqrt{\frac{2 \sigma\E_{\tmop{cap}}
   (| \partial_x h |)}{\rho h}} \tmop{sign} (\partial_x h)$, for the quadratic capillary model. }
    \label{fig:g1dwquad}
\end{figure}

\begin{figure}[htbp]
    \centering
    \begin{subfigure}[t]{0.48\textwidth}
        \centering
        \includegraphics[width=\textwidth]{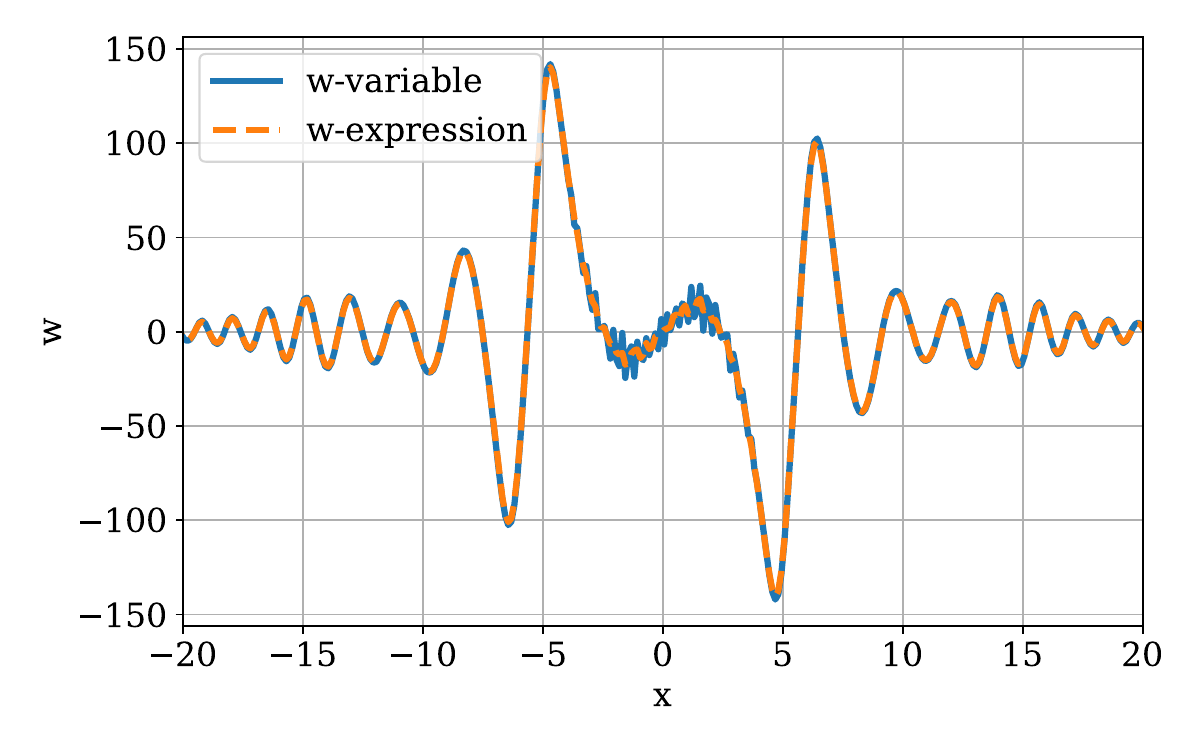}
        \caption{800 cells}
    \end{subfigure}
    \hfill
    \begin{subfigure}[t]{0.48\textwidth}
        \centering
        \includegraphics[width=\textwidth]{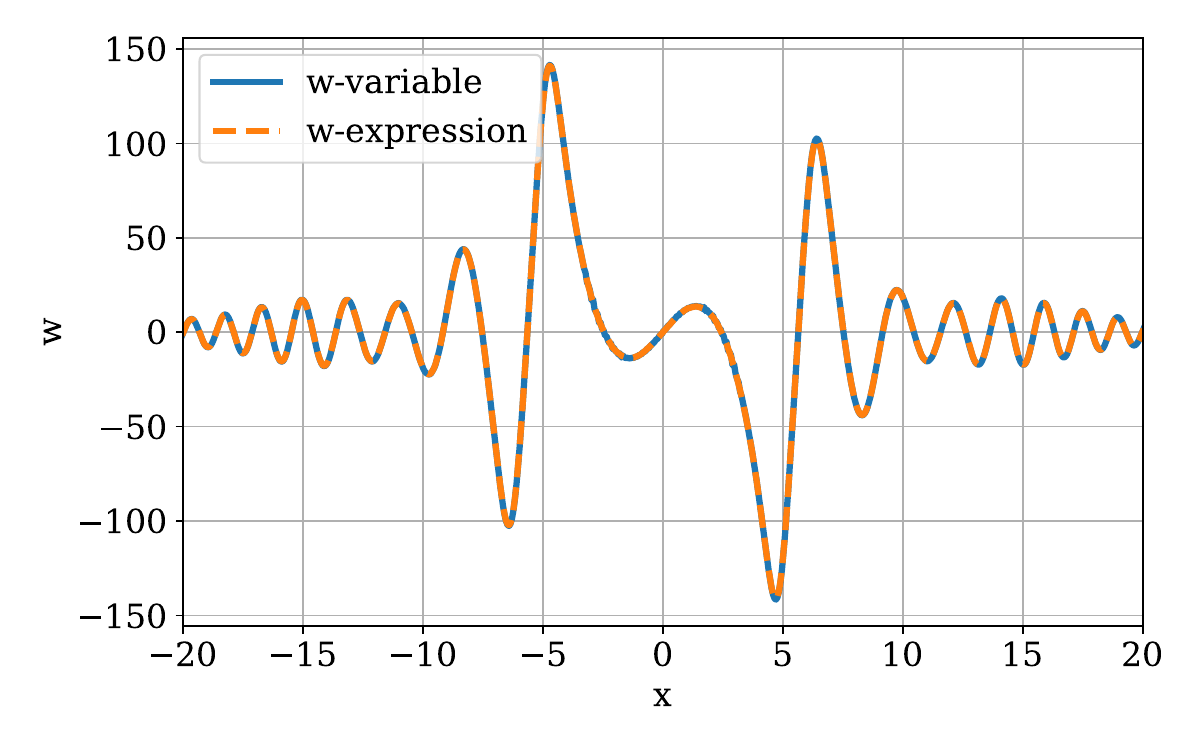}
        \caption{1600 cells}
    \end{subfigure}
    \caption{1-D Gaussian bump:  comparison of $w$ vs $\sqrt{\frac{2 \sigma\E_{\tmop{cap}}
   (| \partial_x h |)}{\rho h}} \tmop{sign} (\partial_x h)$, for the nonlinear capillary model. }
    \label{fig:g1dwnon}
\end{figure}

The quantity $w$ is related to the gradient of the height field but we evolve it as an independent variable. Figure~\ref{fig:g1dwquad} and~\ref{fig:g1dwnon} show the comparison of the computed $w$ and the definition in terms of $h_x$ at the final time for the quadratic and nonlinear cases, respectively. There is good agreement between the two values, which improves under grid refinement. The nonlinear case requires finer grids, as seen in the previous graphs, to avoid the spurious oscillations. Table~\ref{tab:g1dwerr} shows the norm of the difference between these two quantities which decreases at a second-order rate under grid refinement.

\begin{table}[h]
  \centering
  \begin{tabular}{cccccc}
    \toprule
    \multirow{2}{*}{\#cells}
    & \multirow{2}{*}{$\Delta x$}
    & \multicolumn{2}{c}{quadratic}
    & \multicolumn{2}{c}{nonlinear} \\
    \cmidrule(lr){3-4}\cmidrule(lr){5-6}
    & & $L^1$ & $L^\infty$ & $L^1$ & $L^\infty$ \\
    \midrule
      100 & $1.00\times10^{-3}$ & $5.06\times10^{-4}$ & $9.45\times10^{-2}$ & $5.97\times10^{-4}$ & $7.91\times10^{-2}$ \\
      200 & $5.00\times10^{-4}$ & $2.19\times10^{-4}$ & $2.88\times10^{-2}$ & $2.66\times10^{-4}$ & $3.40\times10^{-2}$ \\
      400 & $2.50\times10^{-4}$ & $6.73\times10^{-5}$ & $5.85\times10^{-3}$ & $1.05\times10^{-4}$ & $2.23\times10^{-2}$ \\
      800 & $1.25\times10^{-4}$ & $1.91\times10^{-5}$ & $1.42\times10^{-3}$ & $5.16\times10^{-5}$ & $1.39\times10^{-2}$ \\
     1600 & $6.25\times10^{-5}$ & $4.95\times10^{-6}$ & $3.52\times10^{-4}$ & $9.82\times10^{-6}$ & $1.53\times10^{-3}$ \\
     3200 & $3.13\times10^{-5}$ & $1.27\times10^{-6}$ & $8.87\times10^{-5}$ & $2.32\times10^{-6}$ & $1.38\times10^{-4}$ \\
    \bottomrule
  \end{tabular}
 \caption{1-D Gaussian bump: error between the computed $w$ and its definition as $\sqrt{\frac{2 \sigma\E_{\tmop{cap}} (| \partial_x h |)}{\rho h}} \tmop{sign} (\partial_x h)$}
  \label{tab:g1dwerr}
\end{table}
All the numerical solutions, irrespective of grid size, preserve the total energy in the domain up to machine precision.
\subsection{Gaussian bump in 2-D}

\begin{figure}
	\centering
	\begin{subfigure}[t]{0.48\textwidth}
		\centering
		\includegraphics[width=\textwidth]{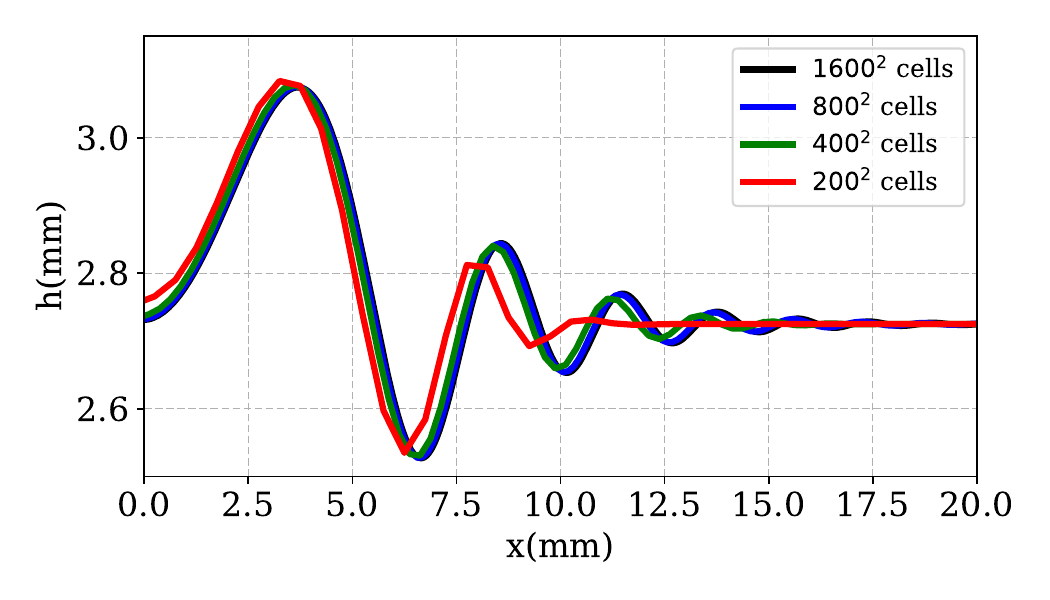}
		\caption{Quadratic}
	\end{subfigure}
	\hfill
	\begin{subfigure}[t]{0.48\textwidth}
		\centering
		\includegraphics[width=\textwidth]{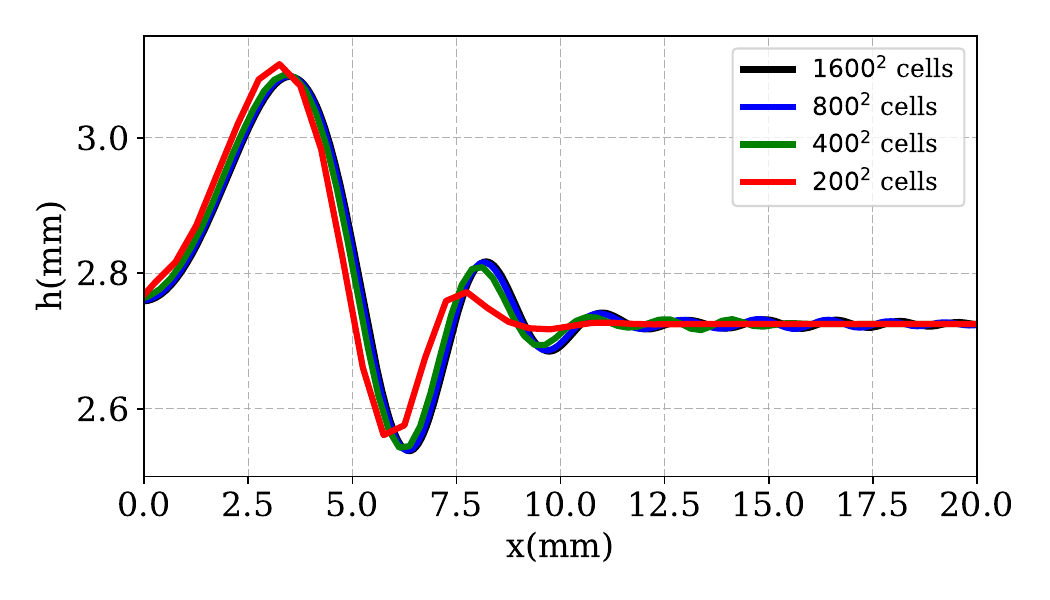}
		\caption{Nonlinear}
	\end{subfigure}\\
	\begin{subfigure}[t]{0.48\textwidth}
		\centering
		\includegraphics[width=\textwidth]{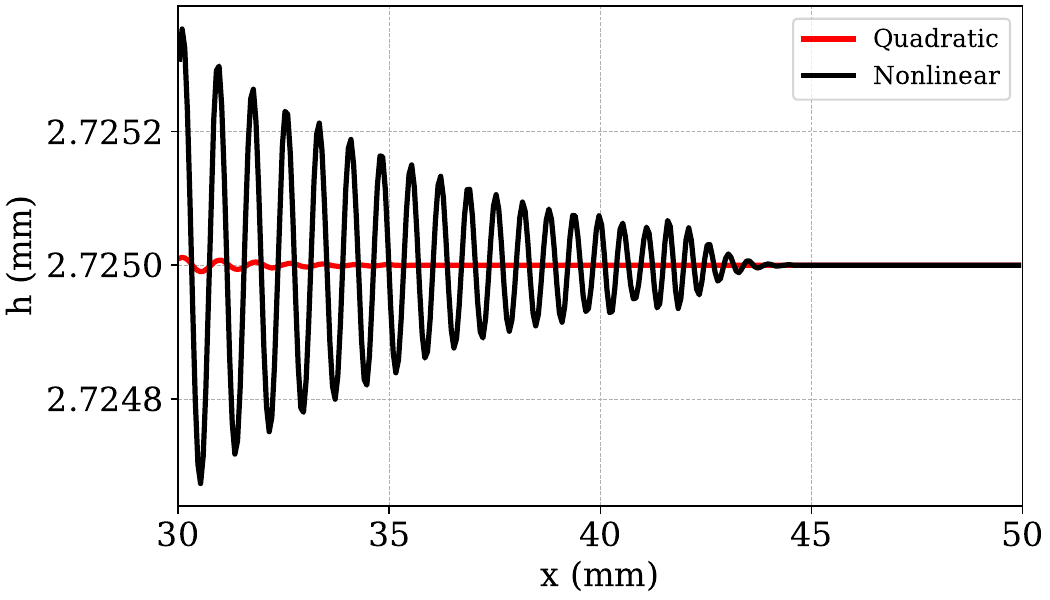}
		\caption{$1600 \times 1600$ cells, $x \in [30 \tmop{mm}, 50 \tmop{mm}]$}
	\end{subfigure}
	\caption{2-D Gaussian bump: solution along the radial line $y=0$.}
	\label{fig:g2dnx}
\end{figure}

\begin{figure}[htbp]
	\centering
	\begin{subfigure}[t]{0.8\textwidth}
		\centering
		\includegraphics[width=\textwidth]{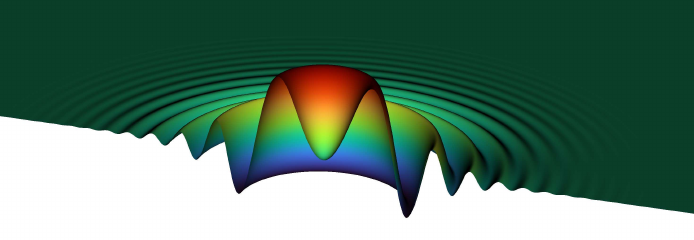}
		\caption{Quadratic capillary}
	\end{subfigure}\\
	\begin{subfigure}[t]{0.8\textwidth}
		\centering
		\includegraphics[width=\textwidth]{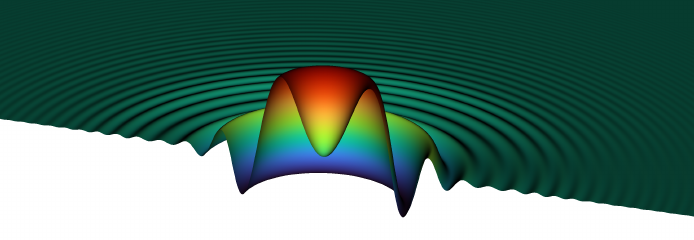}
		\caption{Nonlinear capillary}
	\end{subfigure}
	\caption{ 2-D Gaussian bump: solution computed on $1600\times 1600$ cells, for (a) the quadratic capillary model, (b) the nonlinear capillary model.}
	\label{fig:g2dsurf}
\end{figure}

This is a two-dimensional version of the one-dimensional test case from Section~\ref{sec:gauss1d}. The domain is $[- 50 \tmop{mm}, 50 \tmop{mm}] \times [- 50 \tmop{mm}, 50 \tmop{mm}]$ and the initial condition is
\[
h (x, y, 0) = h_0 + h_1 \exp \left( - \frac{x^2 + y^2}{2 (b / b_0)^2}
\right), \qquad v_1 (x, y, 0) = v_2 (x, y, 0) = 0
\]
and the other parameters are identical to those used in Section~\ref{sec:gauss1d}. The initial condition depends only on the radial coordinate, and thus this is a truly two-dimensional problem. Under the time evolution, the Gaussian mound collapses due to gravity and waves propagate radially outward. The solution along a radial line at the final time is shown in Figure~\ref{fig:g2dnx} for different grid resolutions; the solutions converge under grid refinement. Figure~\ref{fig:g2dnx}(c) shows the zoomed view of the free surface for $x\in [30 \tmop{mm}, 50 \tmop{mm}]$. Similar to the one-dimensional case, we observe a larger amplitude of the waves for the nonlinear case. Figure~\ref{fig:g2dsurf} shows the plot of the free surface at the final time for the quadratic and nonlinear capillary models on a grid of $1600 \times 1600$ cells.

All the numerical solutions, irrespective of grid size, preserve the total energy in the domain up to machine precision.
\subsection{Liquid film falling down an inclined plane: Liu--Gollub experiment}
\label{sec:liugollub}

This is based on an experiment by Liu--Gollub~\cite{Liu1994a,Kalliadasis2012} which consists of a thin film
of an aqueous solution of glycerin on a two-meter-long channel which is inclined at
an angle of $\theta = 6.4^\circ$. The flow is initially uniform and unperturbed, and it is
perturbed periodically at the left end of the domain. Let $h_N, v_N, L, L /
v_N$ be the scales for $h, v, x, t$, respectively, where
\[ v_N = \frac{g h_N^2 \sin \theta}{3 \nu} \]
is the Nusselt speed. The relevant non-dimensional parameters are
\[ \tmop{Re} = \frac{h_N v_N}{\nu}, \qquad \tmop{Fr} = \frac{v_N}{\sqrt{g
   h_N}}, \qquad \tmop{We} = \frac{\rho h_N v_N^2}{\sigma} \]
that characterize the problem. The physical parameters are taken as
\[ \sigma = 0.067~\tmop{N}/\tmop{m}, \qquad \rho = 1080~\tmop{kg}/\tmop{m}^3, \qquad g = 9.80~\tmop{m}/\tmop{s}^2, \qquad \nu = 6.28 \times 10^{- 6}~\tmop{m}^2/\tmop{s} \]
The domain is $[0, 2~\tmop{m}]$ and the initial condition is
\[ h (x, 0) = h_N, \qquad v (x, 0) = v_N \]
where $h_N = 1.279 \tmop{mm}$. At the left boundary $x = 0$, we apply a
periodic forcing in the height field at a frequency of $f$ Hz,
\[
h (0, t) = h_N [1 + 0.1 \sin (2 \pi f t)], \qquad v (0, t) = v_N
\]
and Neumann conditions for all variables on the right boundary. At the left
boundary, the values in ghost cells are set by extrapolation
\[ h_{- 1} = 2 h (0, t) - h_0, \qquad v_{- 1} = 2 v (0, t) - v_0, \qquad w_{-
   1} = 2 w (0, t) - w_0 \]
where
\[ w (0, t) = w \left( h (0, t), \frac{h_0 - h (0, t)}{\Delta x / 2} \right)
\]
The amplitude of the waves exiting from the right boundary is large, and the use of a Neumann condition is not correct. This leads to wave reflections that corrupt the solution inside the domain. We use a sponge layer near the exit where the waves will be damped so that we can use Neumann boundary conditions. As shown in Appendix~\ref{sec:lgbcout}, the use of Neumann conditions does not inject energy into the system at the outflow boundary.

\paragraph{Sponge layer.} Consider an extended domain $[0, 2 + L]$ with a
sponge layer in the region $2 \leqslant x \leqslant 2 + L$. We add a
relaxation term
\[ \frac{\U^{n + 1}_j - \U_j^n}{\Delta t} + \cdots = \cdots + \frac{1}{\tau}
   \eta (x_j) (\U_{\tmop{in}} - \U_j^{\nph}), \qquad \U_{\tmop{in}} =
   \left[\begin{array}{c}
     h_N\\
     h_N v_N\\
     0
   \end{array}\right] \]
and
\[ \eta (x) = \left\{\begin{array}{ll}
     0, & 0 \leqslant x \leqslant 2\\
     \frac{1}{L} (x - 2), & 2 \leqslant x \leqslant 2 + L
   \end{array}\right., \qquad \tau = C_s \frac{L}{v_N} \]
This function is zero at the beginning of sponge layer and takes the value one at the right boundary. The effect of the relaxation term is to drive the solution towards
$\U_{\tmop{in}}$ as the flow moves down the inclined plane. The waves get
damped so that it is easy to apply a Neumann-type outflow boundary condition at
$x = 2 + L$. We use $L = 1$ and $C_s = 0.01$ in the computations.

\begin{figure}[htbp]
    \centering
    \begin{subfigure}[t]{1\textwidth}
        \centering
        \includegraphics[width=\textwidth]{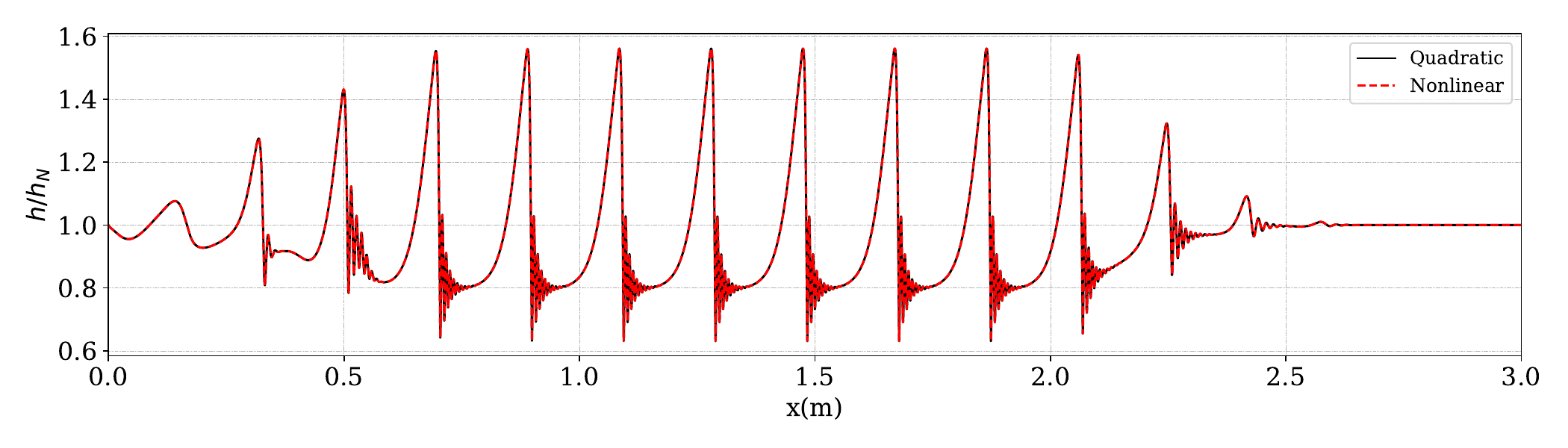}
        \caption{20,000 cells, sponge layer for $x\in[2,3]$}
    \end{subfigure}\\
    \begin{subfigure}[t]{0.48\textwidth}
        \centering
        \includegraphics[width=\textwidth]{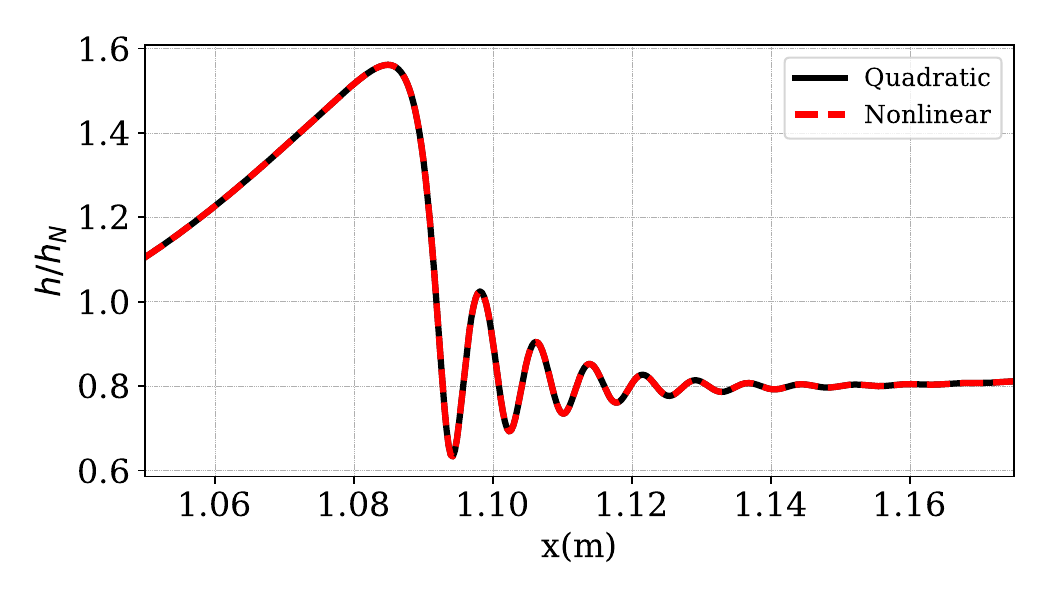}
        \caption{10,000 cells}
    \end{subfigure}
    \hfill
    \begin{subfigure}[t]{0.48\textwidth}
        \centering
        \includegraphics[width=\textwidth]{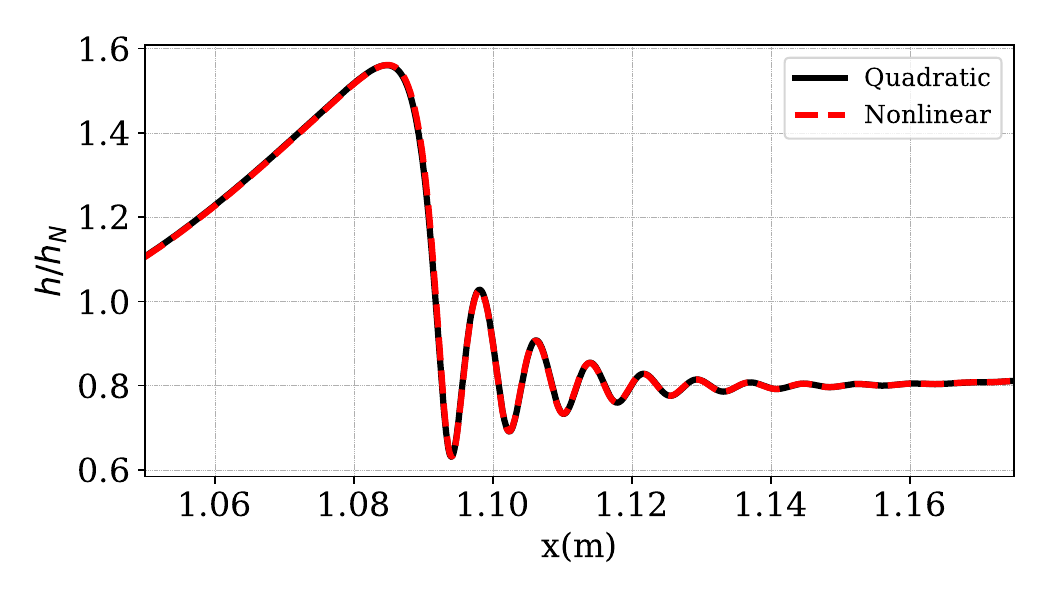}
        \caption{20,000 cells}
    \end{subfigure}
    \caption{Liu--Gollub test case: comparison of the quadratic and nonlinear capillary models for $f=1.5$ Hz.}
    \label{fig:lg1da}
\end{figure}

Figure~\ref{fig:lg1da} shows the solution for frequency $f=1.5$ Hz at the final time in the whole domain; there is little difference between the quadratic and nonlinear capillary models. We also see that the waves are damped in the region $x \in [2,3]$ due to the presence of the sponge layer, and the solution is constant at the right boundary $x=3$ where the use of Neumann conditions does not cause any spurious wave reflections. The zoomed view of the solutions on one of the waves shows the primary wave and the dispersive waves; the solution is well resolved in terms of capturing all the waves even with 10,000 cells.

\begin{figure}[htbp]
    \centering
    \begin{subfigure}[t]{0.48\textwidth}
        \centering
        \includegraphics[width=\textwidth]{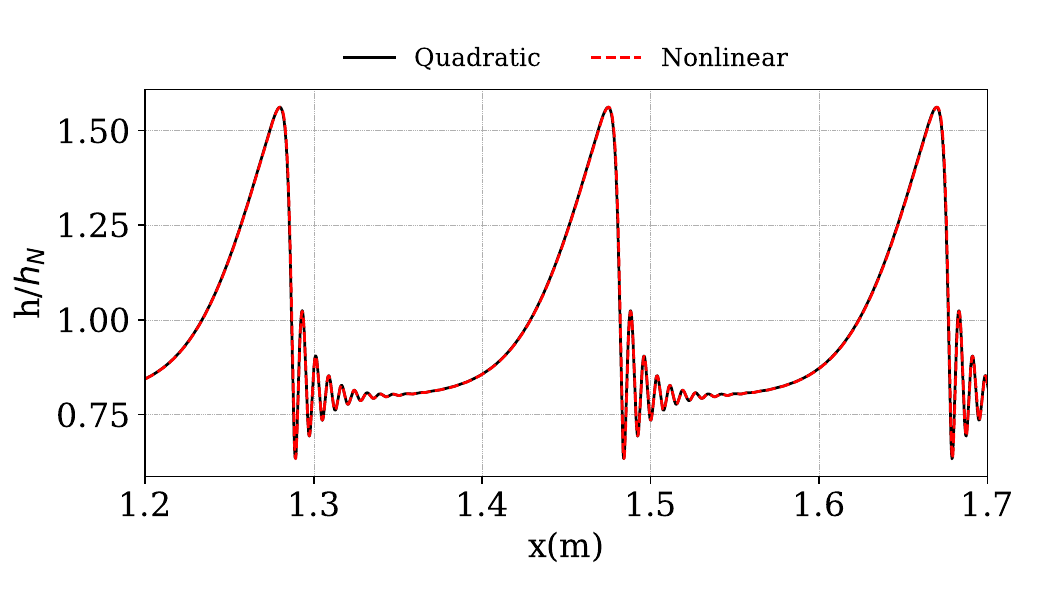}
        \caption{$f=1.5$ Hz}
    \end{subfigure}
    \hfill
    \begin{subfigure}[t]{0.48\textwidth}
        \centering
        \includegraphics[width=\textwidth]{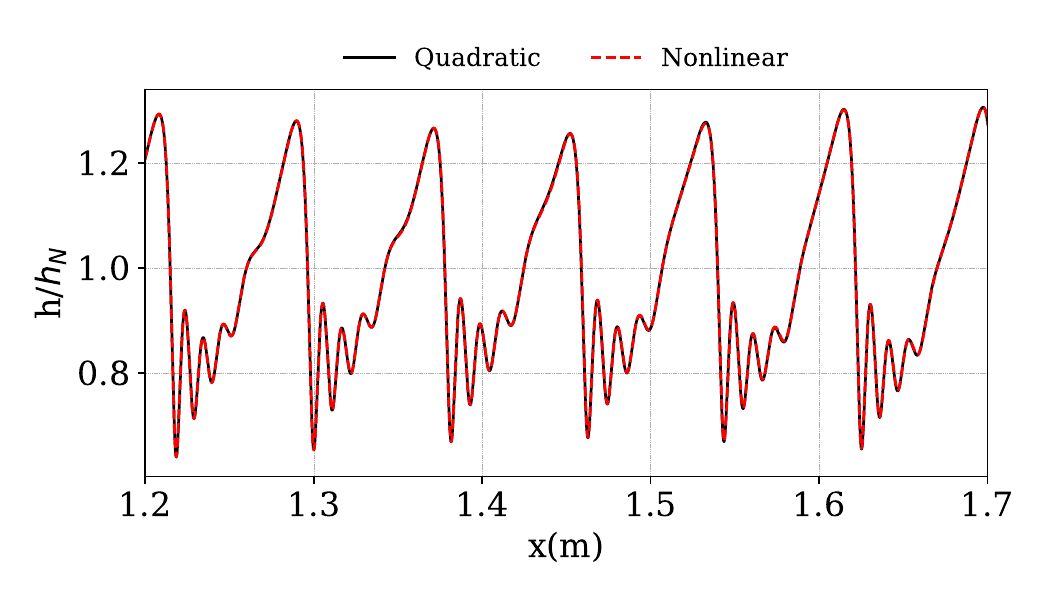}
        \caption{$f=3$ Hz}
    \end{subfigure}\\
        \begin{subfigure}[t]{0.48\textwidth}
        \centering
        \includegraphics[width=\textwidth]{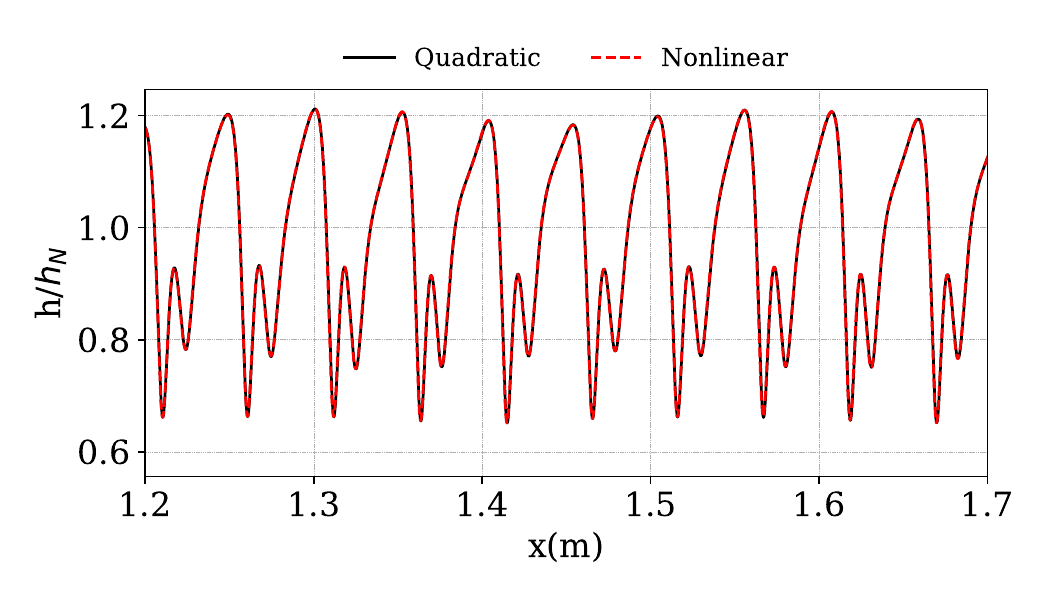}
        \caption{$f=4.5$ Hz}
    \end{subfigure}
    \caption{Liu--Gollub test case: solutions with different frequencies, computed on $10{,}000$ cells.}
    \label{fig:lg1dfreq}
\end{figure}

The qualitative nature of the waves depends on the forcing frequency as observed in the experiments of Liu--Gollub. Figure~\ref{fig:lg1dfreq} shows the solution with three different forcing frequencies, $f=1.5, 3, 4.5$ Hz. As the frequency is increased, the dispersive waves become less prominent and this behavior is observed in the experiments; see Figure~4 in~\cite{Dhaouadi2022} for experimental results.

\begin{figure}
    \centering
    \includegraphics[width=0.9\linewidth]{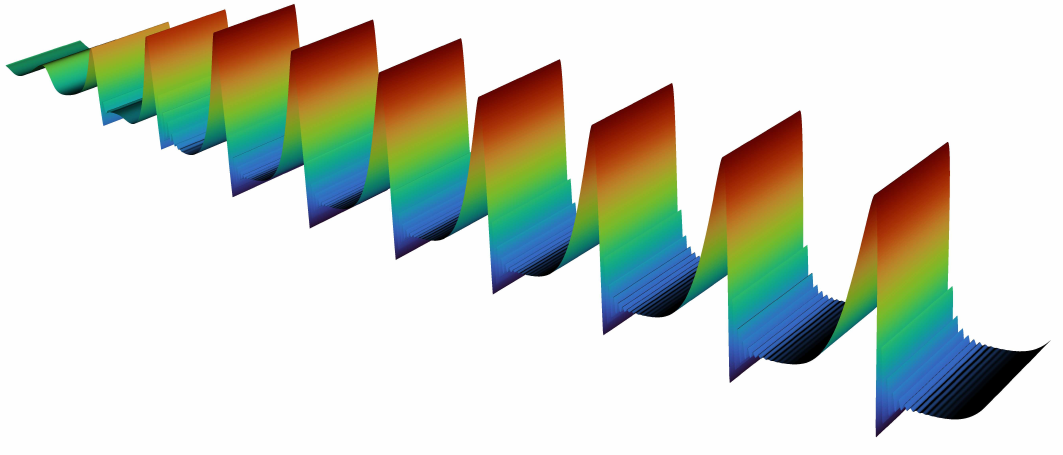}
    \caption{2-D Liu--Gollub experiment: solution in the domain
    $[0,3]\times[-0.25,0.25]$ at time $t = 20$, computed on a grid of $4200 \times 700$ cells and the quadratic capillary model}
    \label{fig:lg2dsurf}
\end{figure}

\paragraph{2-D case.} We repeat this test case in two dimensions on a domain of $[0, 2 + L] \times [- 0.25, 0.25]$, with $L=1$, and the same boundary conditions at $x = 0$ and $x = 2 + L$, and periodic boundary conditions in the $y$ direction. The same idea of a sponge layer is used for $x \in [2, 2 + L]$. Figure~\ref{fig:lg2dsurf} shows the free surface at the final time for the quadratic capillary model, where the flow is from left to right. We observe the same type of waves as seen in the one-dimensional computations.

\subsection{Liquid film falling down a vertical plate: Traveling wave
solutions}

\begin{figure}[htbp]
\centering
\begin{subfigure}[t]{1\textwidth}
\centering
\includegraphics[width=\textwidth]{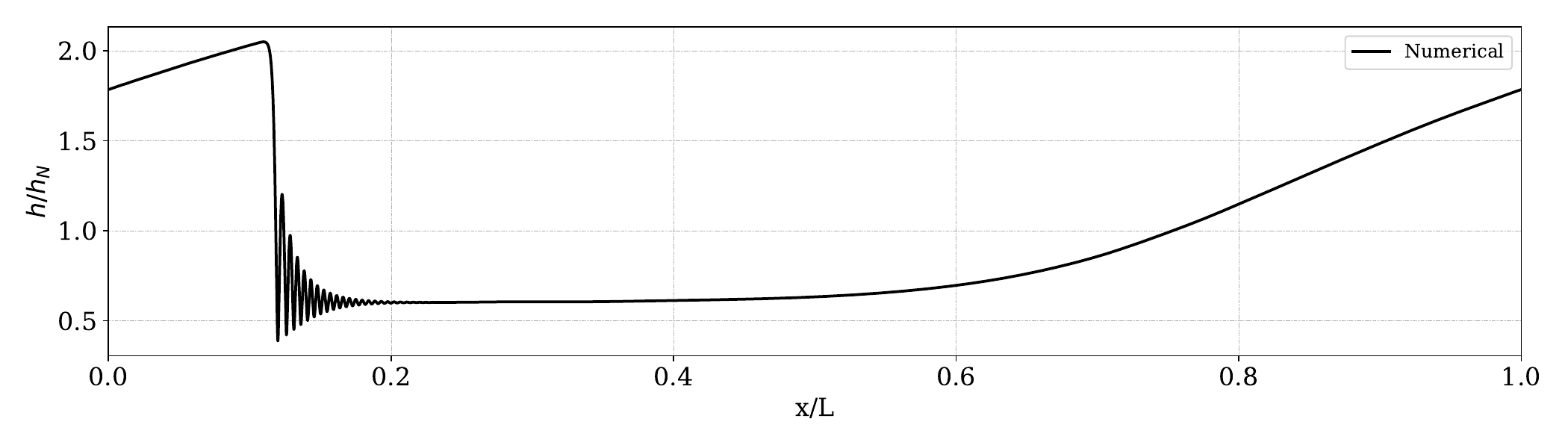}
\caption{Quadratic}
\end{subfigure}\\
\begin{subfigure}[t]{1\textwidth}
\centering
\includegraphics[width=\textwidth]{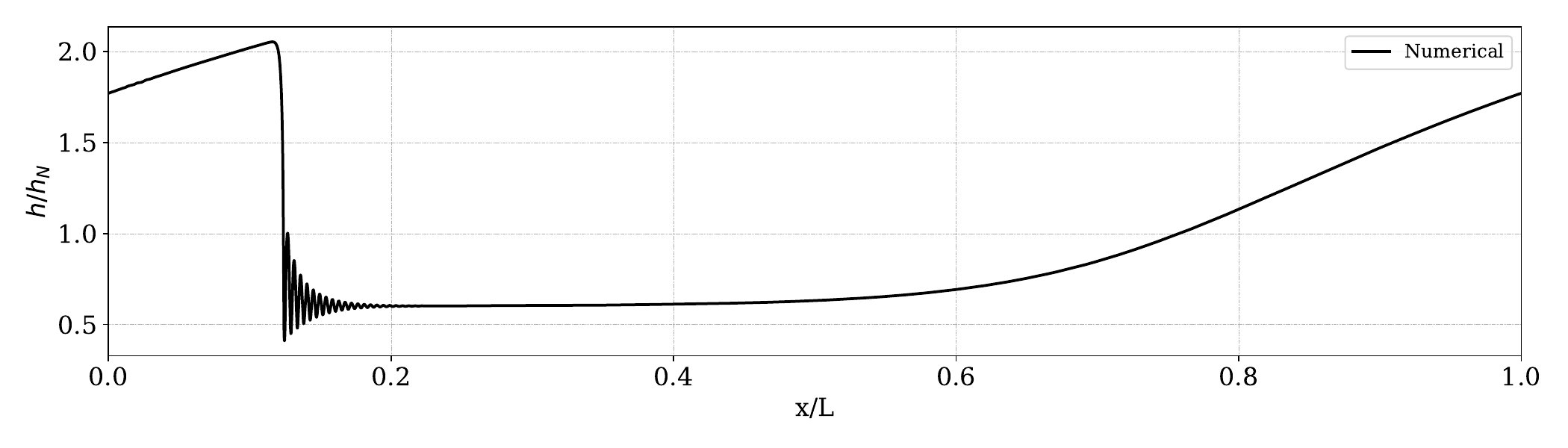}
\caption{Nonlinear}
\end{subfigure}
\caption{Liquid film falling down a vertical plate: solution on 10,000 cells at $t=1~\tmop{s}$}
\label{fig:vfull}
\end{figure}

\begin{figure}
\begin{subfigure}[t]{0.48\textwidth}
\centering
\includegraphics[width=\textwidth]{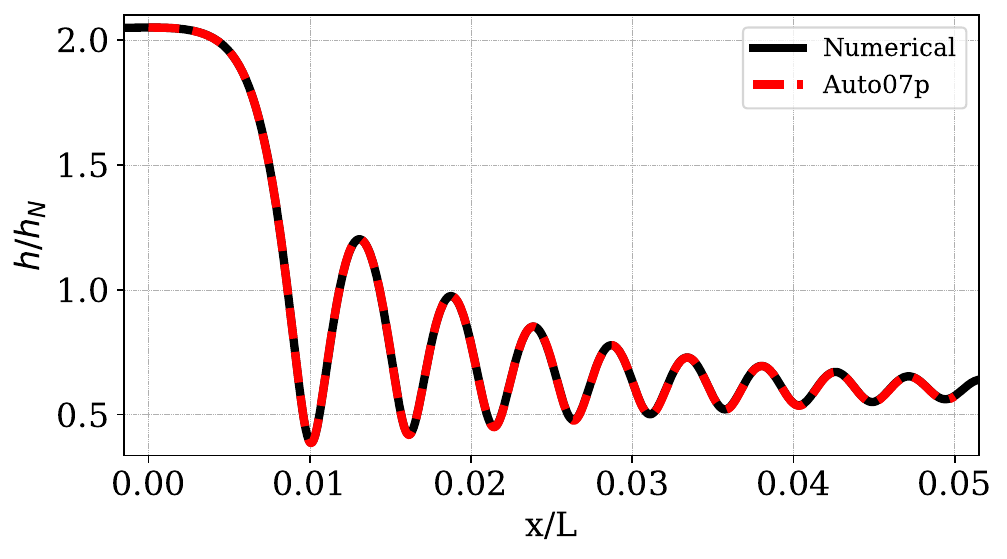}
\caption{Quadratic}
\end{subfigure}\hfill
\begin{subfigure}[t]{0.48\textwidth}
\centering
\includegraphics[width=\textwidth]{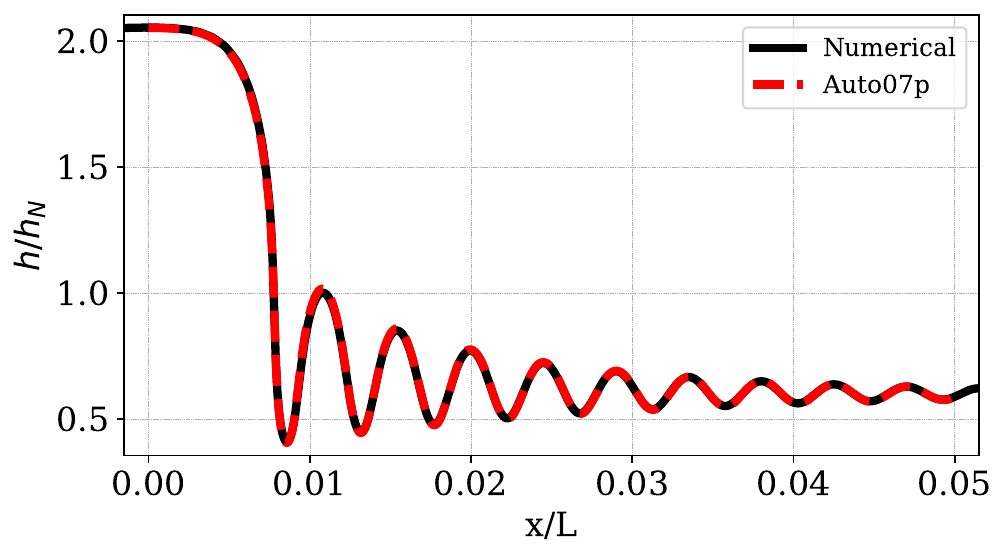}
\caption{Nonlinear}
\end{subfigure}
\caption{Liquid film falling down a vertical plate: zoomed view of the solution compared with AUTO-07P}
\label{fig:vzoom}
\end{figure}

This test case was used in~\cite{Bresch2020} to study the ability of the
model to develop traveling wave solutions. The physical situation involves a
thin liquid film falling down a vertical plate, and we use periodic boundary
conditions. The domain is $[0, L]$ with $L = 400 h_N = 110.636 \tmop{mm}$ and
the initial condition has a sinusoidal perturbation of the height field
\[ h (x, 0) = h_N [1 + 0.1 \sin (2 \pi x / L)], \qquad v (x, 0) = v_N \]
where $h_N = 0.27659 \tmop{mm}$. The other parameters are
\[
\sigma = 0.019322~\tmop{N}/\tmop{m}, \qquad \rho = 994.3~\tmop{kg}/\tmop{m}^3, \qquad g = 9.81~\tmop{m}/\tmop{s}^2, \qquad \nu = 0.931 \times 10^{- 6}~\tmop{m}^2/\tmop{s}
\]
From this initial condition, we evolve the solution up to $t=10~\tmop{s}$. The solution develops a hydraulic jump modulated by dispersive waves and, at large times, assumes an invariant shape that propagates at constant speed. Figure~\ref{fig:vfull} shows the solution at time $t=1~\tmop{s}$, which consists of a single solitary wave with a main hump followed by a train of dispersive waves. Figure~\ref{fig:vzoom} shows a zoomed view of the height profile together with the reference solitary-wave solution computed with AUTO-07P~\cite{AUTO07P} and reported in~\cite{Bresch2020}. In that work, the two-equation model was reduced to an autonomous ODE system in the moving frame $\xi = x-ct$, and solitary waves were obtained by pseudo-arclength continuation of limit cycles from a Hopf bifurcation of the uniform film, under a constant mean-thickness constraint. The numerical and reference profiles are in very good agreement for both capillary models.

\section{Summary and conclusions}\label{sec:conclusion}
Building on the augmented second-order formulation of~\cite{Bresch2020}, we have developed an energy-conserving, second-order, implicit method for both one- and two-dimensional shallow-water models with capillary effects. The main ingredients are an energy-conserving two-point numerical flux derived via Tadmor's
entropy-stability theory, summation-by-parts finite differences, and a special Crank-Nicolson-type scheme, which allows us to prove energy conservation for the fully discrete scheme in the absence of dissipative effects. We consider both quadratic and nonlinear capillary models in a unified method. A wide range of test cases shows that this leads to a stable method that can capture the dispersive waves. Since we use a central-type flux, there is no inherent numerical dissipation; this can lead to spurious numerical oscillations on coarse grids where steep gradients may not be well resolved, especially for the nonlinear capillary model. We can add numerical dissipation to the numerical flux, which yields smooth solutions on coarse grids, but these solutions will not be accurate and do not capture dispersive waves. For smooth solutions, as in all the problems considered here, it is better to use a fine enough grid to eliminate spurious oscillations and capture the dispersive waves, rather than adding numerical dissipation.

This work opens several promising avenues for future research. In particular, adapting this energy-consistent formulation to unstructured grids would greatly broaden its application to complex topographies. Furthermore, incorporating strong magnetic-field coupling would make this augmented framework an invaluable asset for accurate MHD simulation of liquid-metal blankets in nuclear fusion reactors.
\appendix
\section{Results on summation-by-parts}
We state and prove some SBP results that are used in the main part of the paper.
\subsection{Results in 1-D}
For any two grid function $a,b$, we have\footnote{This is the discrete analogue of $\int [a \partial_x b + b \partial_x a] = \int \partial_x (a b) = 0$, assuming periodicity.}
\begin{eqnarray}
  \sum_j [a_j \delta \mu b_j + b_j \delta \mu a_j] & = & \sum_j \left[ a_j
  \frac{b_{j + 1} - b_{j - 1}}{2 \Delta x} + b_j \frac{a_{j + 1} - a_{j -
  1}}{2 \Delta x} \right] \nonumber\\
  & = & \frac{1}{2 \Delta x} \sum_j [a_j b_{j + 1} - a_j b_{j - 1} + b_j a_{j
  + 1} - b_j a_{j - 1}] \nonumber\\
  & = & \frac{1}{2 \Delta x} \sum_j [a_j b_{j + 1} - a_{j + 1} b_j + b_j a_{j
  + 1} - b_{j + 1} a_j] \nonumber\\
  &=&  0  \label{eq:sbp1}
\end{eqnarray}
and\footnote{This is the discrete analogue of $\int a \partial_x [h \partial_x
b] \mathd x = - \int h (\partial_x a) (\partial_x b) \mathd x$, assuming
periodicity.}
\begin{eqnarray*}
  \sum_j a_j \delta [(\mu h) \delta b]_j & = & \frac{1}{\Delta x} \sum_j [a_j
  (\mu h)_{\jph} (\delta b)_{\jph} - a_j (\mu h)_{\jmh} (\delta
  b)_{\jmh}]\\
  & = & \frac{1}{\Delta x} \sum_j [a_j (\mu h)_{\jph} (\delta b)_{\jph} - a_{j + 1} (\mu h)_{\jph} (\delta b)_{\jph}]\\
  & = & - \sum_j (\mu h)_{\jph} (\delta a)_{\jph} (\delta b)_{\jph}
\end{eqnarray*}
The quantity on the right is symmetric with respect to interchanging $a$ and $b$; thus we conclude that
\begin{equation}
  \sum_j a_j \delta [(\mu h) \delta b]_j = \sum_j b_j \delta [(\mu h) \delta
  a]_j \label{eq:sbp2}
\end{equation}

\subsection{Results in 2-D}
For any two grid functions $a,b$, we have
\begin{equation}\label{eq:sbp2da}
(\delta_1 \mu_1 (h \delta_2 \mu_2 a), b)_\Delta = (\delta_2 \mu_2 (h \delta_1
   \mu_1 b), a)_\Delta 
\end{equation}
This can be shown as follows. Starting from the left hand side,
\begin{eqnarray*}
  (\delta_1 \mu_1 (h \delta_2 \mu_2 a), b)_\Delta & = & \frac{\Delta x \Delta y}{2
  \Delta x} \sum_{j, k} \left[ h_{j + 1, k} \frac{a_{j + 1, k + 1} - a_{j + 1,
  k - 1}}{2 \Delta y} - h_{j - 1, k} \frac{a_{j - 1, k + 1} - a_{j - 1, k -
  1}}{2 \Delta y} \right] b_{j, k}\\
  & = & \frac{1}{4} \sum_{j, k} [h_{j + 1, k} a_{j + 1, k + 1} b_{j, k} -
  h_{j + 1, k} a_{j + 1, k - 1} b_{j, k} \nobracket\\
  &  & \hspace{4em} \nobracket - h_{j - 1, k} a_{j - 1, k + 1} b_{j, k} +
  h_{j - 1, k} a_{j - 1, k - 1} b_{j, k}]\\
  & = & \frac{1}{4} \sum_{j, k} [h_{j, k - 1} a_{j, k} b_{j - 1, k - 1} -
  h_{j, k + 1} a_{j, k} b_{j - 1, k + 1} \nobracket\\
  &  & \hspace{4em} \nobracket - h_{j, k - 1} a_{j, k} b_{j + 1, k - 1} +
  h_{j, k + 1} a_{j, k} b_{j + 1, k + 1}]\\
  & = & \frac{\Delta x \Delta y}{2 \Delta y} \sum_{j, k} \left[ h_{j, k + 1}
  \frac{b_{j + 1, k + 1} - b_{j - 1, k + 1}}{2 \Delta x} - h_{j, k - 1}
  \frac{b_{j + 1, k - 1} - b_{j - 1, k - 1}}{2 \Delta x} \right] a_{j, k}\\
  & = & (\delta_2 \mu_2 (h \delta_1 \mu_1 b), a)_\Delta
\end{eqnarray*}

For a vector field like the velocity $\bm{v}$ and a scalar quantity $h$, we have
\begin{equation}\label{eq:sbpdiv}
(\delta_1 \mu_2 (h \nobracket \tmop{div} \bm{v}) \nobracket, v_1)_\Delta +
   (\delta_2 \mu_1 (h \nobracket \tmop{div} \bm{v}) \nobracket, v_2)_\Delta 
   = - \sum_{j, k} h_{\jph, \kph} (\tmop{div}
   \bm{v})_{\jph, \kph}^2 \Delta x \Delta y 
\end{equation}
where $h_{\jph,\kph}$ and $(\tmop{div}\bm{v})_{\jph,\kph}$ are defined in~\eqref{eq:hdivvert} and lead to
\[
h_{\jph,\kph}
=
\frac14\left(
h_{j,k}+h_{j+1,k}+h_{j,k+1}+h_{j+1,k+1}
\right)
\]
and
\begin{align}
(\tmop{div}\bm{v})_{\jph,\kph}
&=
\frac{
v_{1,j+1,k}+v_{1,j+1,k+1}
-v_{1,j,k}-v_{1,j,k+1}
}{2\Delta x}
\nonumber\\
&\quad+
\frac{
v_{2,j,k+1}+v_{2,j+1,k+1}
-v_{2,j,k}-v_{2,j+1,k}
}{2\Delta y}
\label{eq:divvert}
\end{align}
respectively. Denote the two terms in the left hand side of~\eqref{eq:sbpdiv} by $E_x, E_y$ respectively. Under periodic boundary conditions, we shift the summation indices in $E_x$ and collect the contributions at each corner. This gives
\begin{align*}
E_x
&=
(\delta_1 \mu_2 (h \nobracket \tmop{div} \bm{v}) \nobracket, v_1)_\Delta \\
&=
\frac{\Delta y}{2}
\sum_{j,k}
v_{1,j,k}
\Big[
h_{\jph,\kph}
(\tmop{div}\bm{v})_{\jph,\kph}
+h_{\jph,\kmh}
(\tmop{div}\bm{v})_{\jph,\kmh}
\nonumber\\
&\hspace{2.4cm}
-h_{\jmh,\kph}
(\tmop{div}\bm{v})_{\jmh,\kph}
-h_{\jmh,\kmh}
(\tmop{div}\bm{v})_{\jmh,\kmh}
\Big]\\
&=
\frac{\Delta y}{2}\sum_{j,k}
(h\,\tmop{div}\bm{v})_{\jph,\kph}
\Big[
(v_{1,j,k}-v_{1,j+1,k})
+
(v_{1,j,k+1}-v_{1,j+1,k+1})
\Big]
\nonumber\\
&=
-\sum_{j,k}
(h\,\tmop{div}\bm{v})_{\jph,\kph}
\frac{
v_{1,j+1,k}+v_{1,j+1,k+1}
-v_{1,j,k}-v_{1,j,k+1}
}{2\Delta x}
\,\Delta x\,\Delta y
\end{align*}
Similarly, shifting the summation indices in $E_y$ gives
\begin{align*}
E_y
&=
(\delta_2 \mu_1 (h \nobracket \tmop{div} \bm{v}) \nobracket, v_2)_\Delta \\
&= \frac{\Delta x}{2}
\sum_{j,k}
v_{2,j,k}
\Big[]
h_{\jph,\kph}
(\tmop{div}\bm{v})_{\jph,\kph}
+h_{\jmh,\kph}
(\tmop{div}\bm{v})_{\jmh,\kph}
\nonumber\\
&\hspace{2.4cm}
-h_{\jph,\kmh}
(\tmop{div}\bm{v})_{\jph,\kmh}
-h_{\jmh,\kmh}
(\tmop{div}\bm{v})_{\jmh,\kmh}
\Big]
\\
&=
\frac{\Delta x}{2}\sum_{j,k}
(h\,\tmop{div}\bm{v})_{\jph,\kph}
\Big[
(v_{2,j,k}-v_{2,j,k+1})
+
(v_{2,j+1,k}-v_{2,j+1,k+1})
\Big]
\nonumber\\
&=
-\sum_{j,k}
(h\,\tmop{div}\bm{v})_{\jph,\kph}
\frac{
v_{2,j,k+1}+v_{2,j+1,k+1}
-v_{2,j,k}-v_{2,j+1,k}
}{2\Delta y}
\,\Delta x\,\Delta y
\end{align*}
Adding the last two results and using~\eqref{eq:divvert} proves~\eqref{eq:sbpdiv}.
\section{Outflow boundary condition for Liu--Gollub test case}
\label{sec:lgbcout}
The Liu--Gollub~\cite{Liu1994a} test case considered in Section~\ref{sec:liugollub} has inflow and outflow boundaries. While energy is injected into the system through the forcing at the inflow boundary, there should be no energy injection at the outflow boundary. To examine the effect of the Neumann condition used at the outflow boundary, we look at the contribution
of the capillary terms to the total energy in one-dimension, which is given by
\[ \Delta x \sum_{j = 0}^{N - 1} \V_j \cdot \M_j = \Delta x \sum_j \{ v_j \delta
   [(\mu h) \delta (\TensS w)]_j - w_j \TensS_j \delta [(\mu h) \delta v]_j \} - \Delta
   x \sum_j [v_j \delta \mu (\gs w)_j + w_j \gs_j \delta \mu v_j] \]
Note that for any two grid functions $a,b$, we have
\begin{eqnarray}
  \Delta x \sum_{j = 0}^{N - 1} [a_j \delta \mu b_j + b_j \delta \mu a_j] & =
  & \sum_{j = 0}^{N - 1} \left[ a_j \frac{b_{j + 1} - b_{j - 1}}{2} + b_j
  \frac{a_{j + 1} - a_{j - 1}}{2} \right]\nonumber\\
  & = & \frac{1}{2} \sum_{j = 0}^{N - 1} [a_j b_{j + 1} - a_j b_{j - 1} + b_j
  a_{j + 1} - b_j a_{j - 1}]\nonumber\\
  & = & \frac{1}{2} \sum_{j = 0}^{N - 1} [a_j b_{j + 1} + b_j a_{j + 1}]\nonumber\\
  &  & - \frac{1}{2} \sum_{j = - 1}^{N - 2} [a_{j + 1} b_j + b_{j + 1} a_j]\nonumber\\
  & = & - \frac{1}{2} [a_{j + 1} b_j + b_{j + 1} a_j]_{j = - 1} + \frac{1}{2}
  [a_j b_{j + 1} + b_j a_{j + 1}]_{j = N - 1}\nonumber\\
  & = & - \frac{1}{2} [a_0 b_{- 1} + b_0 a_{- 1}] + \frac{1}{2} [a_{N - 1}
  b_N + a_N b_{N - 1}] \label{eq:sbpab1}
\end{eqnarray}
and
\begin{eqnarray}
  \Delta x \sum_{j = 0}^{N - 1} a_j \delta [(\mu h) \delta b]_j & = & \sum_{j
  = 0}^{N - 1} [a_j (\mu h)_{\jph} (\delta b)_{\jph} - a_j (\mu
  h)_{\jmh} (\delta b)_{\jmh}]\nonumber\\
  & = & \sum_{j = 0}^{N - 1} [a_j (\mu h)_{\jph} (\delta b)_{\jph}]
  - \sum_{j = - 1}^{N - 2} [a_{j + 1} (\mu h)_{\jph} (\delta b)_{\jph}]\nonumber\\
  & = & - \Delta x \sum_{j = 0}^{N - 2} (\mu h)_{\jph} (\delta a)_{\jph} (\delta b)_{\jph}\nonumber\\
  &  & + a_{N - 1} (\mu h)_{\Nmh} (\delta b)_{\Nmh} - a_0 (\mu
  h)_{- \half} (\delta b)_{- \half} \label{eq:sbpab2}
\end{eqnarray}
Applying the two SBP formulae~\eqref{eq:sbpab1},~\eqref{eq:sbpab2}, we get
\begin{eqnarray*}
  \Delta x \sum_{j = 0}^{N - 1} \V_j \cdot \M_j & = & v_{N - 1} (\mu h)_{\Nmh} (\delta \TensS w)_{\Nmh} - v_0 (\mu h)_{- 1 / 2} (\delta \TensS w)_{- 1 / 2}\\
  &  & - (\TensS w)_{N - 1} (\mu h)_{\Nmh} (\delta v)_{\Nmh} + (\TensS w)_0
  (\mu h)_{- 1 / 2} (\delta v)_{- 1 / 2}\\
  &  & + \frac{1}{2} [v_0 \gs_{- 1} w_{- 1} + v_{- 1} \gs_0 w_0] - \frac{1}{2}
  [v_{N - 1} \gs_N w_N + v_N \gs_{N - 1} w_{N - 1}]
\end{eqnarray*}
The flux due to capillary terms through the right boundary is
\[ v_{N - 1} (\mu h)_{\Nmh} (\delta \TensS w)_{\Nmh} - (\TensS w)_{N - 1} (\mu
   h)_{\Nmh} (\delta v)_{\Nmh} - \frac{1}{2} [v_{N - 1} \gs_N w_N +
   v_N \gs_{N - 1} w_{N - 1}] \]
which is consistent with the continuous case. If we use Neumann conditions
\[ (h, v, w)_N = (h, v, w)_{N - 1} \]
then the boundary flux, assuming the quadratic capillary model, is
\[ - \frac{1}{2} v_{N - 1} w_{N - 1} [\gs_N + \gs_{N - 1}] = - \frac{1}{2} h_{N -
   1} w_{N - 1}^2 v_{N - 1} \]
If we have outflow, $v_{N - 1} > 0$, then this term is stable in the sense that it decreases the energy. A similar conclusion holds for the nonlinear capillary model.
\section*{Acknowledgments}

BN acknowledges the support grant from the CNRS research federation FR-FCM and EUROfusion. PC acknowledges support from the Department of Atomic Energy, Government of India, under Project Identification No. 1303/9/2025-R\&D-II-DAE/TIFR-17312 and MATRICS Grant No. ANRF/ARGM/2025/001214/MTR from ANRF.

\bibliographystyle{siam}
\bibliography{references}
\end{document}